\documentclass{article}
\usepackage{bbold,pstricks}
\usepackage{amsopn}

\newtheorem{theorem}{Theorem}
\newtheorem{proposition}{Proposition}
\newtheorem{corollary}{Corollary}
\newtheorem{lemma}{Lemma}
\newtheorem{definition}{Definition}

\def\komment#1{}

\DeclareSymbolFont{Extrasymb}{U}{msa}{m}{n}
\DeclareMathSymbol\square\mathrel{Extrasymb}{"03}
\def\qed{{\leavevmode\unskip\nobreak\hfil\penalty 50\hskip 1em%
  \hbox{}\nobreak\hfil\lower 1pt\hbox{$\square $\kern-.5pt}\parfillskip 0pt
  \finalhyphendemerits 0\par}}

\newcommand{\nc}{\newcommand}
\nc{\ep}{\vspace{0.3cm}}
\nc{\Tr}{\mbox{\bf Tr}}
\nc{\nb}{\nabla}
\nc{\na}{\nabla}
\nc{\si}{\sigma}
\nc{\Si}{\Sigma}
\nc{\sid}{\si_D}

\nc{\LL}{{\cal L}}
\nc{\SiM}{\Sigma M}
\nc{\SiN}{\Sigma N}
\nc{\RSi}{R^\Sigma}
\nc{\nbM}{\nabla^M}
\nc{\nbN}{\nabla^N}
\nc{\End}{\mbox{\rm End}}
\nc{\<}{\left\langle}
\renewcommand{\>}{\right\rangle}
\nc{\Spin}{\mbox{\rm Spin}}
\nc{\SO}{\mbox{\rm SO}}
\nc{\spin}{\mbox{\goth spin}}
\nc{\so}{\mbox{\goth so}}
\nc{\noi}{\noindent}
\nc{\DM}{D^M}
\nc{\DN}{D^N}
\nc{\DB}{\mbox{\goth D}^B}
\nc{\Mu}{\mbox{\goth M}}
\nc{\Nu}{\mbox{\goth N}}
\nc{\bis}{,\ldots,}
\nc{\muhut}{{\hat{\mu}}}
\nc{\nbmuhut}{\nabla^\muhut}
\nc{\vp}{\varphi}
\nc{\del}{\partial}
\nc{\eh}{\frac{1}{2}}
\nc{\E}{\mathcal E}
\nc{\udot}{\dot{\cup}}
\nc{\bigudot}{\operatorname*{\stackrel{\raisebox{0pt}{.}}{\bigcup}}}
\nc{\uij}{_{i,j}}
\nc{\Id}{\mbox{Id}}
\nc{\lm}{{L^2(M)}}
\nc{\epnoi}{\ep\noi}
\nc{\Cu}{C^\infty}
\nc{\punkt}{{\hspace{-0.17cm}\bf .}\hspace{0.2cm}}{}
\nc{\intM}{\stackrel{\circ}{M}}
\nc{\barM}{\bar{M}}
\nc{\dM}{\partial M}
\nc{\Adach}{\hat{\mathcal{A}}}
\nc{\NN}{\mathcal{N}}
\nc{\HH}{\mathcal{H}}
\nc{\DD}{\mathcal{D}}
\nc{\BB}{\mathcal{B}}
\nc{\SSS}{\mathcal{S}}
\nc{\Oo}{\mbox{O}}
\nc{\oo}{\mbox{o}}
\nc{\eps}{\epsilon}
\nc{\Hom}{\mbox{Hom}}
\nc{\gt}{{\tilde{g}}}
\nc{\naM}{\nabla^M}
\nc{\naN}{\nabla^N}
\nc{\zz}{{{\Z}/2{\Z}}}

\newfont{\Bbb}{msbm10}
\newfont{\goth}{eufm10}

\begin{document}

\title{The Dirac Operator on Hyperbolic Manifolds of Finite Volume}
\author{Christian B\"ar}
\date{25.~October 2000}
\maketitle

\begin{abstract}
\noindent
We study the spectrum of the Dirac operator on hyperbolic manifolds of
finite volume.
Depending on the spin structure it is either discrete or the whole
real line.
For link complements in $S^3$ we give a simple criterion in terms of
linking numbers for when essential spectrum can occur.
We compute the accumulation rate of the eigenvalues of a sequence of
closed hyperbolic 2- or 3-manifolds degenerating into a noncompact
hyperbolic manifold of finite volume.
It turns out that in three dimensions there is no clustering at all.


\end{abstract}


\setcounter{section}{-1}
\section{Introduction}

The aim of this paper is to study the spectrum of the Dirac operator
on hyperbolic manifolds with finite volume.
Since the corresponding problems for the Laplace-Beltrami operator
acting on differential forms have already been examined let us first
briefly describe those results.
The first natural thing to do is to look at the spectrum of the model
space, $n$-dimensional hyperbolic space $H^n$.
Donnelly \cite{donnelly80a} computed the spectrum of the Laplace operator 
$\Delta_q$ acting on $q$-forms on $H^n$.
For the point spectrum he obtains
$$
spec_p(\Delta_q) = \left\{
\begin{array}{cl}
\{ 0 \}, & q=n/2 \\
\emptyset, & \mbox{otherwise}
\end{array}
\right.
$$
and for the continuous spectrum
$$
spec_c(\Delta_q) = \left\{
\begin{array}{cl}
{[(n-2q-1)^2 /4,\infty)}, & q \le n/2 \\
{[(n-2q+1)^2 /4,\infty)}, & q \ge n/2 
\end{array}
\right.
$$
The eigenvalue $0$ in the case $n=2q$ occurs with infinite multiplicity.
When we pass to quotients of the hyperbolic space we cannot hope to be
able to explicitly compute the spectrum anymore.
But the essential spectrum which is much more robust than the eigenvalues
may still be controlled.
Indeed, Mazzeo and Phillips \cite{mazzeo-phillips90a} showed that except
for the eigenvalue $0$ the essential spectrum on a noncompact hyperbolic 
manifold of finite volume is the same as that of $H^n$
$$
spec_e(\Delta_q) = \left\{
\begin{array}{cl}
{[(n-2q-1)^2 /4,\infty)}, & q \le n/2 \\
{[(n-2q+1)^2 /4,\infty)}, & q \ge n/2 
\end{array}
\right.
$$
In dimension 2 and 3 one can approximate hyperbolic manifolds of finite
volume by compact ones.
In dimension 2 this is clear from Teichm\"uller theory and it can be done
continuously.
In three dimensions it follows from Thurston's cusp closing theorem
that for any noncompact hyperbolic manifold $M$ of finite volume one can 
find a sequence of compact hyperbolic manifolds, pairwise nonhomotopic, 
which converge in a suitable sense to $M$.
What happens to the spectrum under such a degeneration?

Since the spectrum of closed manifolds is discrete we expect that the
eigenvalues in the range of the essential spectrum of the limit manifold
cumulate.
This is true and the rate of clustering has been determined by Ji and 
Zworski \cite{ji-zworski93a} for surfaces, by Chavel and Dodziuk 
\cite{chavel-dodziuk94a} for $n=3$ and $q=0$, and by Dodziuk and McGowan
\cite{dodziuk-mcgowan95a} for $n=3$ and $q=1$.
By Hodge duality this covers all cases.

It turns out that each cusp of the limit manifold $M$ contributes to the
accumulation rate.
This is not surprising because each cusp contributes to the 
essential spectrum.
Let $M_i$ be the approximating sequence of closed hyperbolic manifolds,
$M_i \to M$.
The cusps of $M$ are approximated by degenerating tubes around short closed
geodesics in $M_i$ of length $\ell\uij \stackrel{i \to \infty}
{\longrightarrow} 0$, $j=1\bis k$, where $k$ is the number of cusps of $M$.
For an operator $L$ on a manifold $N$ and an interval $I \subset \R$ we 
introduce the {\em eigenvalue counting function}
$$
\NN_{L,N} (I) := \sharp (spec(L) \cap I) .
$$
Here eigenvalues have to be counted with multiplicity.
Then the accumulation rate turns out to be
$$
\NN_{\Delta_q,M_i} (I) = c(n,q)\frac{x}{\pi} \sum_{j=1}^k \log(1/\ell\uij) 
+ \Oo_x(1)
$$
where
$$
c(n,q) = \left\{
\begin{array}{cll}
2, 
& n=2, & k=0,2 \\
4, 
& n=2, & k=1\\
1/2, 
& n=3, & k=0,3\\
1, 
& n=3, & k=1,2
\end{array}
\right.
$$
and
$$
I = \left\{
\begin{array}{cll}
{[1/4,1/4+x^2]}, 
& n=2 & \\
{[1,1+x^2]}, 
& n=3, & k=0,3 \\
{[0,x^2]}, 
& n=3, & k=1,2
\end{array}
\right.
$$
and $\Oo_x(1)$ denotes an error term bounded as a function of $i$.
Moreover, Colbois and Courtois \cite{colbois-courtois89a,colbois-courtois91a} 
showed that the eigenvalues below the
bottom of the essential spectrum of $M$ are limits of eigenvalues of the $M_i$.

We want to study the analogous questions for the Dirac operator $D$ acting
on spinors, sometimes also called Atiyah-Singer operator, on hyperbolic
manifolds.
The spectrum of the model space $H^n$ has been computed by Bunke 
\cite{bunke91a}.
Note that there is an incorrect statement about the eigenvalue $0$ in that
paper.
See also \cite{camporesi92a,camporesi-higuchi96a} and the remark after
the proof of Lemma \ref{separiervar} in this paper.
The result is
$$
spec_p(D) = \emptyset,
\hspace{1cm}
spec_c(D) = \R .
$$
Since $D$ is of first order the spectrum is not semibounded.
When we pass to nonsimply connected hyperbolic manifolds a new
piece of structure enters the picture for which there is no analog
for the Laplace operator.
We have to specify a spin structure on the manifold.
First of all, this means that we have to restrict our attention to hyperbolic 
spin manifolds.
In particular, the manifolds must be orientable.
If the manifold is spin the spin structure is not unique.
There are as many spin structures on $M$ as there are elements in the 
cohomology group $H^1(M;\zz)$.
It turns out that the choice of spin structure has dramatic impact on
the Dirac spectrum.
We will define in Section 3 what it means that a spin structure is {\em trivial
along a cusp} of a hyperbolic manifold.
This is an essentially topological property.
Our first result is

{\bf Theorem \ref{dichotomy}.}
{\em Let $M$ be a hyperbolic manifold of finite volume equipped with
a spin structure.

If the spin structure is trivial along at least one cusp, then the Dirac 
spectrum satisfies
$$
spec(D) = spec_e(D) = \R .
$$

If the spin structure is nontrivial along all cusps, then the spectrum
is discrete,
$$
spec(D) = spec_d(D) .
$$
}

We see already that there is no analog for the eigenvalues below the bottom
of the essential spectrum as studied by Colbois and Courtois.

We will see that if $M$ is 2- or 3-dimensional and has only one cusp, then
only the second case occurs, the spectrum is always discrete (Corollary
\ref{einende}).
If $M$ is a surface with at least two cusps, then both cases occur.
The spin structure can be made trivial on any choice of an even number
of cusps.

In three dimensions this is not true in general.
It can happen that the spectrum is always discrete even if the manifold
has more than one cusp.
If the hyperbolic manifold is given as the complement of a link in $S^3$,
then there is simple criterion to decide if there is a spin structure
such that $spec(D) = \R$.

{\bf Theorem \ref{verschlingung}.}
{\em Let $K \subset S^3$ be a link, let $M = S^3 - K$ carry a hyperbolic metric
of finite volume.

If the linking number of all pairs of components $(K_i,K_j)$ of $K$
is even,
$$
Lk(K_i,K_j) \equiv 0 \mbox{ mod } 2,
$$
$i \not= j$, then the spectrum of the Dirac operator on $M$ is discrete for all
spin structures,
$$
spec(D) = spec_d(D).
$$

If there exist two components $K_i$ and $K_j$ of $K$, $i \not= j$, with
odd linking number, then $M$ has a spin structure such that the spectrum of the
Dirac operator satisfies
$$
spec(D) = \R.
$$
}

Determining linking numbers modulo 2 is equivalent to counting overcrossings
modulo 2, hence extremely simple.
See the last section for examples.

Next we study the behavior under the degeneration process in dimension 2
and 3.
Of course, we have to assume that the spin structure on $M$ is, in a 
suitable sense, the limit of the spin structures on the $M_i$.
In two dimensions the result is

{\bf Theorem \ref{ausart2d}.}
{\em
Let $M_i$ be a sequence of closed hyperbolic surfaces converging
to a noncompact hyperbolic surface $M$ of finite volume.
Let each $M_i$ have exactly $k$ tubes with trivial spin structure around 
closed geodesics of length $\ell\uij$ tending to zero.
Hence $M$ has exactly $2k$ cusps along which the spin structure is trivial.
Let $x > 0$.

Then the eigenvalue counting function for the Dirac operator satisfies
for sufficiently small $\ell\uij$:
{\em $$
\NN_{D,{M_i}}(-x,x) = \frac{4x}{\pi}\sum_{j=1}^k \log(1/\ell\uij) 
+ \Oo_{x}(1).
$$}
}

Very recently, Farinelli \cite{farinelli98a} gave an upper bound on
the spectral accumulation of the lower part of the Dirac spectrum of 
hyperbolic 3-manifolds.
However, we will show that in three dimensions there is no clustering at all!

{\bf Theorem \ref{ausart3d}.}
{\em
Let $M_i$ be a sequence of closed hyperbolic 3-manifolds converging
to a noncompact hyperbolic 3-manifold $M$ of finite volume.
Let each $M_i$ have exactly $k$ tubes around 
closed geodesics of length $\ell\uij$ tending to zero.
Hence $M$ has exactly $k$ cusps.
Let $x > 0$.

Then the spin structure is nontrivial along all tubes and
the eigenvalue counting function for the Dirac operator remains bounded:
{\em $$
\NN_{D,{M_i}}(-x,x) = \Oo_{x}(1).
$$}
}

The reason for this fact, at first glance quite surprising, is of topological
nature.
The spin structure on the tubes must be nontrivial because the trivial spin 
structure on the 2-torus is nontrivial in spin cobordism $\Omega^{Spin}_2$.
In other words, the spin structures on hyperbolic 3-manifolds of finite
volume for which $spec(D) = \R$ do not occur as limits of spin structures
on closed hyperbolic 3-manifolds.

We see that the freedom to choose different spin structures leads to new
phenomena in the spectral theory of the Dirac operator on hyperbolic
manifolds for which there is no analog for the Laplace operator.
This also distinguishes the classical Dirac operator acting on spinors
from those twisted Dirac operators on locally symmetric spaces which
have typically been studied in the context of representation theory 
\cite{atiyah-schmid77a,atiyah-schmid79a} and index theory \cite{mueller87a}.

The paper is organized as follows.
In the first section we collect a few facts about hyperbolic manifolds.
The structure of the cusps and tubes is important for our purposes.
A description of the degeneration process in dimension 2 and 3 is given.

In the second section we present some generalities about the $L^2$-spectrum
of self-adjoint elliptic operators.
We give a prove of the so-called decomposition principle which roughly
says that modifying the manifold and the operator in a compact region of the
manifold does not affect the essential spectrum.
This will be extremely useful for us because we can restrict our attention
to the cusps of the hyperbolic manifolds.
There are many versions of this principle in the literature but we found it 
convenient to prove it in a quite general form.
Our version can e.g.\ be applied to the Dirac operator on manifolds with
boundary with suitable boundary conditions.

In Section 3 we prove Theorem \ref{dichotomy}.
We use a separation of variables along the cusps which reduces the problem
to the study of simple Schr\"odinger operators on an interval.

In the forth section we derive a general version of domain monotonicity.
This allows one to estimate eigenvalues by cutting the manifold into
pieces.
This has been used extensively in the spectral geometry of the Laplace 
operator.
Here we need this tool for the Dirac operator.

We are then able to prove Theorem \ref{ausart2d} in Section 5.
It is important that tubes in a hyperbolic surface are warped products
so that the separation of variables can again be applied.

We would like to do the same thing in three dimensions in Section 7 
but we have the problem that tubes are no longer simple warped products.
Therefore we include a general formula in Section 6 which relates the
square of the Dirac operator on a manifold foliated by hypersurfaces
to operators along the leaves and normal derivatives.
This way we can regard the square of the Dirac operator on the tube as 
a Schr\"odinger operator acting on Hilbert space-valued functions on an
interval.
We will then be able to prove Theorem \ref{ausart3d} in Section 7.

In the last section we discuss the different spin structures which 2-
or 3-dimen\-sional hyperbolic manifolds can have.
This is more topological in nature.
We conclude with a few examples of link complements for which essential
spectrum does or does not occur.

{\bf Acknowledgements.}
It is a pleasure to thank W.\ Ballmann, J.\ Dodziuk, U.\ Hamenst\"adt, and
W.\ M\"uller for many fruitful discussions and valuable hints.
This paper was written while the author enjoyed the hospitatility of
SFB 256 at the University of Bonn.


\section{Hyperbolic Manifolds of Finite Volume}

A {\em hyperbolic manifold} is a complete connected Riemannian
manifold of constant sectional curvature -1.
We collect a few well-known facts about such manifolds with special
emphasis on the case of finite volume.
A thorough introduction to the topic is given in \cite{benedetti-petronio92a}.
 
Every hyperbolic manifold $M$ of finite volume can be decomposed disjointly 
into a relatively compact $M_0$ and finitely many cusps $\E_j$,
\begin{equation}
M = M_0 \udot \bigudot_{j=1}^k \E_j
\label{zerlegung}
\end{equation}
where each $\E_j$ is of the form $\E_j = N_j \times [0,\infty)$.
Here $N_j$ denotes a connected compact manifold with a flat metric 
$g_{N_j}$, a {\em Bieberbach manifold}, and $\E_j$ carries the warped 
product metric $g_{\E_j} = e^{-2t}\cdot g_{N_j} + dt^2$.

\begin{center}
\pspicture(1,0)(14,6)

\pscustom[linecolor=white,fillstyle=solid,fillcolor=lightgray]{
  \pscurve(4.95,6)(4.9,5.5)(4.7,5)(4.4,4.5)
  \pscurve(4.4,4.5)(4.6,4.4)(5.05,4.3)(5.7,4.5)
  \pscurve(5.7,4.5)(5.4,5)(5.2,5.5)(5.15,6)
  \pscurve(5.15,6)(4.95,6)}
\pscustom[linecolor=white,fillstyle=solid,fillcolor=lightgray]{
  \pscurve(7,3.6)(8,3.35)(9,3.22)(10,3.16)(11,3.13)
  \pscurve(11,3.13)(11,2.97)
  \pscurve(11,2.96)(10.6,2.96)(10,2.94)(9,2.88)(8,2.75)(7,2.5)
  \pscurve(7,2.5)(6.9,2.6)(6.8,3.05)(7,3.6)}

\pscurve(7,2.5)(8,2.75)(9,2.88)(10,2.94)(11,2.97)
\pscurve(7,3.6)(8,3.35)(9,3.22)(10,3.16)(11,3.13)
\pscurve(4.4,4.5)(4.7,5)(4.9,5.5)(4.95,6)
\pscurve(5.7,4.5)(5.4,5)(5.2,5.5)(5.15,6)
\psecurve(5.4,5)(5.7,4.5)(6.15,3.95)(7,3.6)(8,3.35)
\psecurve(4.7,5)(4.4,4.5)(3,3)(3,1)(4.5,0.7)(5.9,1.2)(7,2.5)(8,2.75)
\pscurve(4.1,2.75)(4,2.5)(4.5,2.1)(5,2.2)(5.2,2.3)
\pscurve(4,2.5)(4.6,2.45)(5,2.2)
\pscurve(7,2.5)(6.8,3.05)(7,3.6)
\pscurve(4.4,4.5)(5.05,4.3)(5.7,4.5)

\rput(5.3,2.9){$M_0$}
\rput(7.3,3.1){\psframebox*[framearc=0.5]{$\E_1$}}
\rput(5.05,4.8){\psframebox*[framearc=0.5]{$\E_2$}}
\rput(6.55,3.1){{$N_1$}}
\rput(5.05,4.05){{$N_2$}}

\endpspicture
\vspace{-1cm}
\end{center}
\begin{center}
\bf Fig.~1
\end{center}

For example, $N_j$ could be a flat torus, as is always the case if $M$ is
2- or 3-dimensional and orientable.
This simple structure of the cusps will allow us to apply a separation
of variables technique to the Dirac operator on hyperbolic manifolds
of finite volume.

It turns out that very different phenomena occur in hyperbolic geometry 
depending on the dimension.
In dimension 2 there is a whole continuum of hyperbolic structures
(hyperbolic metrics modulo isometries) on a given surface.
This is known as Teichm\"uller theory.
In particular, if we fix a compact surface $M$, then there are 
continuous deformations of hyperbolic metrics on $M$ under which $M$
degenerates to a noncompact hyperbolic surface of finite volume.
These deformations correspond to paths in the Teichm\"uller space
converging to the boundary.

In contrast, in dimension $n \ge 3$, we know by Mostow's rigidity theorem 
that any compact manifold carries at most one hyperbolic structure.
Therefore continuous degenerations are not possible.

If $n = 3$ however, the following kind of degeneration still occurs.
Thurston's cusp closing theorem says that for any hyperbolic manifold 
$M = M_0 \udot \bigudot_{j=1}^k \E_j$ of finite volume with metric $g$
there are compact hyperbolic manifolds $(M_i,g_i)$ which can 
be decomposed disjointly into
$$
M_i = M_{i,0} \udot \bigudot_{j=1}^k T_{i,j}
$$
where $T_{i,j}$ is the closed tubular neighborhood of radius $R_{i,j}$ 
of a simple closed geodesic $\gamma_{i,j} \subset M_i$ of length $\ell_{i,j}$.
The boundary $N_{i,j} = \del T_{i,j}$ is a flat torus.
In the degeneration $(i \to \infty)$ the following happens:
\begin{itemize}
\item
$\ell_{i,j} \to 0$

\item
$R_{i,j} = \eh\log(1/\ell\uij) + c_0 \to \infty$ where $c_0$ is some constant
independent of $i$.

\item
There are diffeomorphisms $\Phi_i: \bar{M}_{0} \to \bar{M}_{i,0}$ of compact 
manifolds with boundary such that the metrics 
$\Phi_i^\ast(g_i|_{\bar{M}_{i,0}})$ converge in the $C^\infty$-topology
to $g|_{\bar{M}_{0}}$.

\item
The pull-backs of the metrics of $N_{i,j}$ converge in the 
$C^\infty$-topology to the one of $N_j$.

\end{itemize}
Moreover, if we write 
$$
T\uij [r_1,r_2] = \{ x \in M_i\ |\ dist(x,\gamma\uij) \in [r_1,r_2]\ \}
$$ 
for the tubular region around $\gamma\uij$, so that $T\uij = T\uij [0,R\uij]$,
then we have in addition
\begin{itemize}
\item
For every $0 < r_1 < r_2 < R\uij$ the tubular region $T\uij [r_1,r_2]$ is 
isometric to $T^2 \times [r_1,r_2]$ with the metric
$g_{r} + dr^2$ where $g_r$ is the flat metric on the 2-torus given
by the lattice $\Gamma_r \subset \R^2$ spanned by the vectors
$(2\pi\sinh(r),0)$ and $(\alpha_{i,j}\sinh(r), \ell_{i,j}cosh(r))$
for some ``holonomy angle'' $\alpha_{i,j} \in [-\pi,\pi]$.
\end{itemize}

\begin{center}
\pspicture(1,0)(14,6)

\psline[linestyle=dotted](2,1)(12,1)
\psline[linestyle=dotted](4,0)(4,5)
\psline{->}(4,1)(8,1)
\psline{->}(4,1)(7,4)
\psline[linestyle=dotted](4,4)(7,4)
\psline[linestyle=dotted](8,1)(11,4)
\psline[linestyle=dotted](7,4)(11,4)
\psdots(11,4)
\rput(6,0.5){$2\pi\sinh(r)$}
\rput(5.5,4.5){$\alpha\uij\sinh(r)$}
\rput(3,2.5){$\ell\uij\cosh(r)$}

\endpspicture
\end{center}
\begin{center}
\bf Fig.~2
\end{center}

This description of the degeneration is also valid in the 2-dimensional
case, except that the tube $T\uij$ is of the form
$T\uij = S^1 \times [-R\uij,R\uij]$ with metric $ds^2 = \ell\uij^2
\cosh(t)^2 d\theta^2 + dt^2$ where $t \in [-R\uij,R\uij]$, $\theta
\in S^1 = \R/\Z$ and $R\uij = \log(1/\ell\uij) + c_0$.
In particular, the boundary of the
tube is of the form $S^0 \times S^1 = S^1 \udot S^1$.
Hence each tube degenerates into {\em two} cusps.

\begin{center}
\pspicture(1,0)(14,6)

\pscurve(3,5.7)(5,5.3)(7,5.1)(9,5.3)(11,5.7)
\pscurve(3,4.3)(5,4.7)(7,4.9)(9,4.7)(11,4.3)
\psellipse(7,5)(0.05,0.1)
\rput(7,5.4){$\gamma_{i,j}$}
\pscurve(3,5.7)(2.8,5)(3,4.3)
\psecurve[linestyle=dashed](1,6.5)(2,6)(3,5.7)(5,5.3)
\psecurve[linestyle=dashed](1,3.5)(2,4)(3,4.3)(5,4.7)
\psline(2.8,5)(6.95,5)
\rput(3.5,5){\psframebox*[framearc=0.5]{$R_{i,j}$}}
\rput(4,5.7){\psframebox*[framearc=0.5]{$T_{i,j}$}}

\pscurve(11,5.7)(11.2,5)(11,4.3)
\psecurve[linestyle=dashed](13,6.5)(12,6)(11,5.7)(9,5.3)
\psecurve[linestyle=dashed](13,3.5)(12,4)(11,4.3)(9,4.7)

\pscurve(3,2.7)(5,2.3)(7,2.1)(9,2.05)
\pscurve(3,1.3)(5,1.7)(7,1.9)(9,1.95)
\pscurve(3,2.7)(2.8,2)(3,1.3)
\psecurve[linestyle=dashed](1,3.5)(2,3)(3,2.7)(5,2.3)
\psecurve[linestyle=dashed](1,0.5)(2,1)(3,1.3)(5,1.7)

\pscurve(5,1.05)(7,1.1)(9,1.3)(11,1.7)
\pscurve(5,0.95)(7,0.9)(9,0.7)(11,0.3)
\pscurve(11,1.7)(11.2,1)(11,0.3)
\psecurve[linestyle=dashed](13,2.5)(12,2)(11,1.7)(9,1.3)
\psecurve[linestyle=dashed](13,-0.5)(12,0)(11,0.3)(9,0.7)

\rput(7.2,3.5){$\downarrow i \to \infty$}

\endpspicture
\vspace{-1cm}
\end{center}
\begin{center}
\bf Fig.~3
\end{center}

In order to define the Dirac operator we also need to specify spin
structures on our manifolds.
In the degeneration we require that the diffeomorphisms 
$\Phi_i: \bar{M}_{0} \to \bar{M}_{i,0}$ can be chosen compatible with 
the spin structures.


\section{Generalities about the $L^2$-Spectrum}

Let $H$ be a complex Hilbert space and let $A$ be a 
self-adjoint linear operator with dense domain
$A : \mathcal{D}(A) \subset H \to H$.

\begin{definition}\punkt
A number $\lambda \in \C$ is called an {\bf eigenvalue} of $A$ if
$A - \lambda \Id$ is not injective.
In this case we call $\dim(\ker(A - \lambda \Id))$ the {\bf multiplicity}
of $\lambda$.
The set of eigenvalues, $spec_p(A)$, is called the {\bf point spectrum}.

The {\bf essential spectrum}, $spec_e(A)$, is the set of $\lambda \in \C$ 
for which there exists a sequence $x_i \in \mathcal{D}(A)$ satisfying
\begin{equation}
\| x_i \| = 1, 
\hspace{1cm}
(A-\lambda \Id)x_i \to 0,
\hspace{1cm}
x_i \rightharpoonup 0
\label{weyl}
\end{equation}
for $i \to \infty$.
Here ``$\rightharpoonup$'' denotes weak convergence as opposed to
norm convergence ``$\to$''.

The union of the point spectrum and the essential spectrum is the
{\bf spectrum} of $A$, $spec(A) = spec_p(A) \cup spec_e(A)$.
\end{definition}

Note that the spectrum of a self-adjoint operator is actually contained
in $\R$ and that the point spectrum and the essential spectrum need not
be disjoint.
Eigenvalues of infinite multiplicity and eigenvalues
which are cumulation points of the spectrum are contained in both
the point spectrum and the essential spectrum.

\begin{definition}\punkt
The set 
$$
spec_d(A) = spec_p(A) - spec_e(A)
$$ 
is called the {\bf discrete spectrum}.
The set 
$$
spec_c(A) = spec_e(A) - spec_p(A)
$$ 
is called the {\bf continuous spectrum}.
\end{definition}

Sometimes it will be convenient to look at the square of an operator instead
of the operator itself.
We will then use that $spec_e(A) = \emptyset$ if and only if
$spec_e(A^2) = \emptyset$.

In the definition of the essential spectrum (\ref{weyl}) can be replaced
by other equivalent conditions.
For example, instead of demanding $x_i \rightharpoonup 0$ we could require
that there is no convergent subsequence.
If the operator $A$ is the closure of an operator $L$ with domain
$\mathcal{D}(L)$, then we can as well require $x_i \in \mathcal{D}(L)$.
See e.g.\ \cite{weidmann80a} for details.
A sequence as in (\ref{weyl}) is called a {\em Weyl sequence}.

Let us show that the essential $L^2$-spectrum of self-adjoint elliptic 
differential operators on manifolds does not change when one
modifies the manifold in a compact region.

In what follows we will denote the space of $L^p$-sections in
a Hermitian vector bundle $E$ over a Riemannian manifold $M$ by
$L^p(M,E)$, the Sobolev space of sections whose covariant derivatives
up to order $k$ are $L^p$ by $H^{k,p}(M,E)$.
The space of $k$ times continuously differentiable sections is
denoted by $C^k(M,E)$, $0 \le k \le \infty$, and the space of
$C^k$-sections with compact support is denoted by $C^k_0(M,E)$.

\begin{proposition}\punkt
{\em (Decomposition Principle)}
\newline
Let $\bar{M}$ be a Riemannian manifold, with (possibly empty) compact boundary,
$\bar{M} = \intM \udot \del M$.
Let $E$ be a Hermitian vector bundle over $\bar{M}$.
Let $L$ be an essentially self-adjoint linear differential operator 
of order $d \ge 1$ with domain $\mathcal{D}(L)$, 
$C^\infty_0(\intM,E) \subset \mathcal{D}(L) \subset C^\infty_0(\barM,E)$.
Suppose for every compact $K \subset \barM$ there is an {\em elliptic estimate}
\begin{equation}
\|x \|_{H^{d,2}(K,E)} \le
C\cdot \left( \|x \|_{L^2(\bar{M},E)} + \|Lx \|_{L^2(\bar{M},E)}  \right)
\label{ellipticestimate}
\end{equation}
for all $x \in \mathcal{D}(L)$, $C=C(K)$.
Denote the closure of $L$ in $L^2(\bar{M},E)$ by $\bar{L}$.

Let $\bar{M}'$ be another Riemannian manifold and let $E'$, $L'$, 
and $\bar{L'}$ be defined similarly on $\bar{M}'$.
We assume there exist compact sets $K \subset \bar{M}$, $K' \subset \bar{M}'$
such that $\bar{M}-K = \bar{M}'-K'$, and $E=E'$, $L=L'$ over $\bar{M}-K$.

Then 
$$
spec_e(\bar{L}) = spec_e(\bar{L'}).
$$
\label{decompprinciple}
\end{proposition}

\begin{center}
\pspicture(1,0)(14,6)

\pscustom[linecolor=white,fillstyle=solid,fillcolor=lightgray]{
  \pscurve(2.4,4.5)(2.7,5)(2.9,5.5)(2.95,6)
  \pscurve(3.15,6)(3.18,5.8)(3.2,5.5)(3.4,5)(3.7,4.5)
  \pscurve(3.7,4.46)(3.5,4.4)(3.05,4.3)(2.4,4.5)}

\pscustom[linecolor=white,fillstyle=solid,fillcolor=lightgray]{
  \pscurve(5,2.5)(6,2.75)(7,2.88)
  \pscurve(7,3.22)(6.7,3.25)(6,3.35)(5,3.6)
  \pscurve(4.99,3.6)(4.9,3.4)(4.8,3.05)(5,2.5)}

\pscurve(5,2.5)(6,2.75)(7,2.88)
\pscurve(5,3.6)(6,3.35)(7,3.22)
\pscurve(2.4,4.5)(2.7,5)(2.9,5.5)(2.95,6)
\pscurve(3.7,4.5)(3.4,5)(3.2,5.5)(3.15,6)
\psecurve(3.4,5)(3.7,4.5)(4.15,3.95)(5,3.6)(6,3.35)
\psecurve(2.7,5)(2.4,4.5)(1,3)(1,1)(2.5,0.7)(3.9,1.2)(5,2.5)(6,2.75)
\pscurve(2.1,2.75)(2,2.5)(2.5,2.1)(3,2.2)(3.2,2.3)
\pscurve(2,2.5)(2.6,2.45)(3,2.2)
\pscurve(5,2.5)(4.8,3.05)(5,3.6)
\pscurve(2.4,4.5)(3.05,4.3)(3.7,4.5)

\rput(3.2,3.3){$K$}
\rput(4.5,4.4){$M$}


\pscustom[linecolor=white,fillstyle=solid,fillcolor=lightgray]{
  \pscurve(8.4,4.5)(8.7,5)(8.9,5.5)(8.95,6)
  \pscurve(9.15,6)(9.18,5.8)(9.2,5.5)(9.4,5)(9.7,4.5)
  \pscurve(9.7,4.46)(9.5,4.4)(9.05,4.3)(8.4,4.5)}

\pscustom[linecolor=white,fillstyle=solid,fillcolor=lightgray]{
  \pscurve(11,2.5)(12,2.75)(13,2.88)
  \pscurve(13,3.22)(12.7,3.25)(12,3.35)(11,3.6)
  \pscurve(10.99,3.6)(10.9,3.4)(10.8,3.05)(11,2.5)}

\pscurve(11,2.5)(12,2.75)(13,2.88)
\pscurve(11,3.6)(12,3.35)(13,3.22)
\pscurve(8.4,4.5)(8.7,5)(8.9,5.5)(8.95,6)
\pscurve(9.7,4.5)(9.4,5)(9.2,5.5)(9.15,6)
\psecurve(9.4,5)(9.7,4.5)(10.15,3.95)(11,3.6)(12,3.35)
\pscurve(11,2.5)(10.8,3.05)(11,3.6)
\pscurve(8.4,4.5)(9.05,4.3)(9.7,4.5)
\psecurve(9.3,5)(8.4,4.5)(8.5,1.5)(11,2.5)(12,2.75)

\rput(9.6,1.6){\psframebox*[framearc=0.5]{}}
\psccurve(9.5,0.9)(9.85,0.9)(10.1,1.4)(9.6,1.15)
\pscurve(9.48,0.92)(9.55,1.4)(9.3,1.9)
\pscurve(10.1,1.4)(9.7,1.7)(9.6,2)

\rput(9.3,3){$K'$}
\rput(10.5,4.4){$M'$}
\rput(10.4,0.9){$\partial M'$}

\endpspicture
\end{center}
\begin{center}
\bf Fig.~4
\end{center}

Note that sections in $C^\infty_0(\bar{M},E)$ need not vanish on $\dM$
in contrast to those of $C^\infty_0(\intM,E)$.

In case $\dM = \emptyset$ the elliptic estimate (\ref{ellipticestimate})
holds automatically if $L$ is an elliptic operator \cite[p.\ 379, Thm.\ 11.1]
{taylor96a}.
In this case the decomposition principle can be found in many places in the
literature for various operators (mostly of second order), see e.g.\ 
\cite{glazman65a,donnelly-li79a,eichhorn88a}.

In the presence of boundary establishing (\ref{ellipticestimate}) is subtler.
It usually follows from {\em coercive estimates}
$$
\|x \|_{H^{d,2}(K,E)}^2 \le
C\cdot \left( \|x \|_{L^2(\bar{M},E)}^2 + \|Lx \|_{L^2(\bar{M},E)}^2 +
\sum_j \| B_jx\|_{H^{d-d_j-1/2,2}(\dM,E)}^2 \right)
$$
where $B_j$ are boundary (pseudo-) differential operators of order 
$d_j \le d-1$, $x \in \Cu_0(\barM,E)$.
If $B_jx|_{\dM}=0$ for all $x\in\mathcal{D}(L)$, then 
(\ref{ellipticestimate}) holds.
The coercive estimate is automatic if $L$ together with the $B_j$
form a regular elliptic boundary value problem \cite[V.11]{taylor96a}.
For example, a Laplace type operator $L$ together with Dirichlet boundary
conditions $x|_{\dM}=0$ forms a regular elliptic boundary value problem.
We will use Proposition \ref{decompprinciple} with Dirichlet boundary 
conditions for the square of the Dirac operator which by the Lichnerowicz 
formula \cite{lichnerowicz63a} 
\begin{equation}
D^2 = \nb^\ast\nb + \frac{scal}{4}
\label{lichnerowicz}
\end{equation}
is of Laplace type.

Here a {\em Laplace type operator} is an operator of the form
$$
L = \nb^\ast\nb + \Re
$$
where $\nb$ is a metric connection on a Hermitian vector
bundle over a Riemannian manifold, $\nb^\ast$ is its $L^2$-adjoint and
$\Re$ is a smooth symmetric endomorphism field (zero order term).
Laplace type operators are special elliptic operators of second order.

One can also apply the decomposition principle directly to the Dirac operator
with suitable boundary conditions.
Since we will not use this fact we leave the details to the reader.

{\sc Proof of Proposition \ref{decompprinciple}.}
Since the whole situation is symmetric in $\bar{M}$ and $\bar{M}'$ it is 
sufficient to show $spec_e(\bar{L}) \subset spec_e(\bar{L'})$.
Let $\lambda \in spec_e(\bar{L})$ and let $x_i \in \DD(L) \subset
L^2(\bar{M},E)$ be a Weyl sequence as in (\ref{weyl}).

Choose a compact subset $K_1 \subset \bar{M}$ whose interior contains $K$
and another compact subset $K_2 \subset \bar{M}$ whose interior contains $K_1$.
By the elliptic estimate (\ref{ellipticestimate})
\begin{eqnarray*}
\|x_i\|_{H^{d,2}(K_2,E)} 
&\le&
C\cdot \left( \|x_i\|_{L^2(\bar{M},E)} + \|Lx_i\|_{L^2(\bar{M},E)}  \right) \\
&\le&
C \cdot \left( 1 + \|Lx_i - \lambda x_i \|_{L^2(\bar{M},E)}
+ |\lambda | \right) \\
&\stackrel{(\ref{weyl})}{\le}& C'.
\end{eqnarray*}
Since $(x_i)_i$ is bounded in the $H^{d,2}$-norm and $K_2$ is compact we can, 
by the Rellich's lemma, pass to a subsequence, again denoted $(x_i)_i$, 
which converges in
$H^{d-1,2}(K_2,E)$ to some element $x_\infty \in H^{d-1,2}(K_2,E)$.

To compute $x_\infty$ we pick a cut-off function $\psi_1$ identical to 1 
on $K_1$ and vanishing outside $K_2$.
On the one hand, since $\psi_1 x_i \to \psi_1 x_\infty$ in $L^2(K_2,E)$,
$$
(\psi_1 x_i, \psi_1 x_\infty )_{L^2(\bar{M},E)} =
(\psi_1 x_i, \psi_1 x_\infty )_{L^2(K_2,E)} \longrightarrow
(\psi_1 x_\infty, \psi_1 x_\infty )_{L^2(K_2,E)}.
$$
On the other hand,
$$
(\psi_1 x_i, \psi_1 x_\infty )_{L^2(\bar{M},E)} =
(x_i, \psi_1^2 x_\infty )_{L^2(\bar{M},E)} \longrightarrow 0
$$
because $x_i \rightharpoonup 0$.
Hence $\psi_1 x_\infty = 0$ and $x_\infty|_{K_1} = 0 \in H^{d-1,2}(K_1,E)$.
Therefore
\begin{equation}
\|x_i\|_{H^{d-1,2}(K_1,E)} \longrightarrow 0.
\label{kompaktnull}
\end{equation}
In particular, for $i$ sufficiently large, $\|x_i\|^2_{L^2(K_1,E)} \le \eh$ 
and thus 
\begin{equation}
\|x_i\|^2_{L^2(\bar{M}-K_1,E)} \ge \eh .
\label{vonunten}
\end{equation}

Choose a cut-off function $\psi \in C^\infty(\bar{M},\R)$ with
$\psi = 0$ on $K$ and $\psi = 1$ on $\bar{M}-K_1$,
$0 \le \psi \le 1$ everywhere.
Let us look at the sequence $y_i \in L^2(\bar{M}',E')$ where 
$y_i = \psi\cdot x_i$ on $\bar{M}-K = \bar{M}'-K'$ and $y_i \equiv 0$ on $K'$.
First of all, by (\ref{vonunten}), 
$$
\|y_i\|^2_{L^2(\bar{M}',E')} \ge \|x_i\|^2_{L^2(\bar{M}-K_1,E)} \ge \eh.
$$
Secondly, for any $z \in L^2(\bar{M}',E')$,
$$
(y_i,z)_{L^2(\bar{M}',E')} = (x_i,\psi z)_{L^2(\bar{M},E)} \to 0
$$
by (\ref{weyl}).
Hence $y_i \rightharpoonup 0$.

Thirdly,
$$
L'y_i = L(\psi x_i) = \psi Lx_i + Q x_i
$$
where $Q = [L,\psi]$ is a differential operator of order $d-1$.
Moreover, $Q$ vanishes outside $K_1$ because $\na\psi$ does.
There is a constant $C_2 > 0$ such that
$$
\|Q x_i\|_{L^2(\bar{M}',E')} \le
C_2 \cdot \|x_i\|_{H^{d-1,2}(K_1,E)}.
$$
Therefore (\ref{kompaktnull}) implies $\|Q x_i\|_{L^2(\bar{M}',E')}
\to 0$.
We conclude
\begin{eqnarray*}
\| L'y_i - \lambda y_i \|_{L^2(\bar{M}',E')} 
&\le&
\| \psi (Lx_i - \lambda x_i) \|_{L^2(\bar{M},E)}  + 
\|Q x_i\|_{L^2(\bar{M}',E')} \\
&\le&
\| Lx_i - \lambda x_i \|_{L^2(\bar{M},E)}  + 
\|Q x_i\|_{L^2(\bar{M}',E')} 
\longrightarrow 0.
\end{eqnarray*}
Thus the sequence $(y_i/||y_i||_{L^2(\bar{M'},E')})_i$ is a Weyl sequence for 
the operator $\bar{L'}$.
Hence $\lambda \in spec_e(\bar{L'})$.
\qed

The proposition will be very useful for the study of the essential
spectrum of hyperbolic manifolds because it tells us that we only
need to consider the operator on the cusps and those have a very simple
form.


\section{The Dirac Operator on Hyperbolic Manifolds}

In this section we will study the type of the spectrum of the
Dirac operator on hyperbolic manifolds of finite volume.
Studying the type means finding out if the spectrum is e.g.\ purely
discrete or purely essential or contains both components.
Our Dirac operator will always be the classical Dirac operator,
sometimes also called Atiyah-Singer operator, acting on spinors.
For definitions see \cite{lawson-michelsohn89a}.

If $M$ is an $n$-dimensional Riemannian spin manifold and $N \subset M$
is an oriented hypersurface, then every spin structure on $M$ canonically
induces a spin structure on $N$.
If $n$ is odd, then the restriction to $N$ of the spinor bundle $\SiM$
of $M$ is precisely the spinor bundle of $N$, $\SiM|_N = \SiN$.
If $n$ is even, then $\SiM|_N$ is isomorphic to $\SiN\oplus\SiN$.

Let $H$ denote the mean curvature function of $N$ with respect to
the unit normal field $\nu$.
Let $\DM$ be the Dirac operator of $M$.
Let $\DN$ be the Dirac operator of $N$ in case $n$ is odd.
If $n$ is even let $\DN$ be the direct sum of the Dirac operator
of $N$ and its negative.
In either case $\DN$ acts on sections of $\SiM|_N$.
The two operators $\DM$ and $\DN$ are related by the formula
\begin{equation}
-\nu \cdot D^M \sigma = D^N \sigma - \frac{n-1}{2} H \sigma + \nabla^M_{\nu} 
\sigma ,
\label{diracgauss}
\end{equation}
see e.g.\ \cite{baer96b,trautman95a}.
Here $\sigma$ is a section of $\SiM$ defined in a neighborhood of $N$,
``$\cdot$'' denotes Clifford multiplication with respect to the 
manifold $M$ and $\nbM$ is the Levi-Civita connection of $\SiM$.

The case of a warped product will be of special importance.
Let $N$ be an $(n-1)$-dimensional Riemannian spin manifold, let
$I \subset \R$ be an interval.
We give $M = N \times I$ the product spin structure and the warped
product metric
$$
ds^2(x,t) = \rho(t)^2g_N(x) + dt^2
$$
where $\rho : I \to \R$ is a fixed positive smooth function.
For example, cusps of a hyperbolic manifold are of this form with
$I = [0,\infty)$ and $\rho(t) = e^{-t}$.
Let $\nu = \frac{\del}{\del t}$ be the unit vector field along $I$.
The mean curvature of $N \times \{ t \}$ in $M$ is now given by 
$H(t) = -\frac{\dot{\rho}(t)}{\rho(t)}$.

\begin{lemma}\punkt
\label{separiervar}
Let $M$ be a warped product as above.
Suppose there is a subspace $X$ of the kernel of $D^N$ such that
$$
\ker(D^N) = X \oplus \nu\cdot X, 
\hspace{1cm}
X \perp \nu\cdot X .
$$
Write $d = \dim(X) = \dim(\ker(D^N))/2$.
Let $0 < \mu_1 \le \mu_2 \le \mu_3 \cdots \to \infty$ be the positive 
eigenvalues of $D^N$, each eigenvalue repeated according to its multiplicity.

Then there is a unitary equivalence 
$$
L^2(M,\SiM) \to \bigoplus_{\mu\in spec(D^N)} L^2(I,\C,dt) =
\bigoplus_{j=1}^d L^2(I,\C^2,dt) \oplus
\bigoplus_{j=1}^\infty L^2(I,\C^2,dt)
$$
under which the Dirac operator $D^M$ is transformed into
$$
D^M \to \bigoplus_{j=1}^d D_0 \oplus
\bigoplus_{j=1}^\infty D_{\mu_j}
$$
where
$$
D_\mu = \left(  
\begin{array}{cc}
0 & -\frac{d}{dt} + \frac{\mu}{\rho(t)} \\
\frac{d}{dt} + \frac{\mu}{\rho(t)} & 0
\end{array}
\right) .
$$
Similarly, the square of the Dirac operator is transformed into
$$
(D^M)^2 \to \bigoplus_{\mu\in spec(D^N)} L_\mu
$$
where
$$
L_\mu = -\frac{d^2}{dt^2} + \frac{\mu \dot{\rho}(t)}{\rho(t)^2}
+ \frac{\mu^2}{\rho(t)^2}
$$
on $L^2(I,\C,dt)$.
\end{lemma}

{\sc Proof.}
We decompose
$$
L^2(N,\SiM|_N) =
\HH^+ \oplus X \oplus \nu\cdot X \oplus \HH^- ,
$$
where $\HH^\pm$ is the sum of eigenspaces of $D^{N}$ for
positive or negative eigenvalues respectively.
Let $\vp_1,\vp_2,\vp_3,\ldots$ be orthonormal eigenvectors corresponding
to the positive eigenvalues $0 < \mu_1 \le \mu_2 \le \mu_3 \cdots \to \infty$,
$\vp_j \in L^2(N,\SiM|_N)$.
Then we have the Hilbert space decomposition
$$
\HH^+ = \bigoplus_{j=1}^\infty \C\cdot\vp_j .
$$
Since Clifford multiplication with $\nu$ anticommutes with $D^N$
we see that $\nu\cdot\vp_j$ is an eigenvector for the eigenvalue
$-\mu_j$ and hence
$$
\HH^- = \bigoplus_{j=1}^\infty \C\cdot\nu\cdot\vp_j .
$$
Write $\HH_j := \C\cdot\vp_j \oplus \C\cdot\nu\cdot\vp_j$.
Similarly, let $\psi_1 \bis \psi_d$ be an orthonormal basis of $X$
and put $\tilde{\HH}_j := \C\cdot\psi_j \oplus \C\cdot\nu\cdot\psi_j$.
Then
$$
L^2(N,\SiM|_N) = \bigoplus_{j=1}^d \tilde{\HH}_j\oplus 
\bigoplus_{j=1}^\infty \HH_j .
$$
By (\ref{diracgauss}) the Dirac operator $D^M$ leaves the
Hilbert space decomposition
$$
L^2(M,\SiM) = \bigoplus_{j=1}^d L^2(I,\tilde{\HH}_j,\rho(t)^{n-1}dt)
\oplus \bigoplus_{j=1}^\infty L^2(I,\HH_j,\rho(t)^{n-1}dt)
$$
invariant and $D^M(\alpha_j\vp_j + \alpha_{-j}\nu\vp_j) =
(\alpha_{-j}\mu_j + \frac{n-1}{2}H\alpha_{-j} - \dot{\alpha}_{-j} ) \vp_j +
(\alpha_j\mu_j - \frac{n-1}{2}H\alpha_j + \dot{\alpha}_j ) \nu\vp_j$.
The map 
$$
\alpha_j\vp_j + \alpha_{-j}\nu\vp_j \mapsto
\rho(t)^{\frac{n-1}{2}}\left(
\begin{array}{c}
\alpha_j \\
\alpha_{-j}
\end{array}
\right)
$$
yields a unitary equivalence
$$
L^2(I,\HH_j,\rho(t)^{n-1}dt) \to L^2(I,\C^2,dt)
$$
under which the Dirac operator is transformed into
$$
D_{\mu_j} =
\left(
\begin{array}{cc}
0 & -\frac{d}{dt} + \frac{\mu_j}{\rho(t)}  \\
\frac{d}{dt} + \frac{\mu_j}{\rho(t)} & 0
\end{array}
\right)
$$
and similarly for the zero eigenvalues.
The formula for the square of the Dirac operator follows immediately.
\qed

{\bf Remark.}
The assumption $\ker(D^N) = X \oplus \nu\cdot X$, $X \perp \nu \cdot X$,
is necessary only for the decomposition of the Dirac operator $D^M$
itself, not for its square $(D^M)^2$.
This assumption is automatically satisfied if $M$ has even dimension.
In this case $\SiM|_N = \Si^+M|_N \oplus \Si^-M|_N \cong \Si N \oplus
\Si N$ and one can simply take $X = \ker(D^N) \cap \Cu(N,\Si^+M|_N)$.
If $\dim(M)$ is odd, then the assumption is equivalent to $\Adach(N) = 0$.

{\bf Remark.}
Lemma \ref{separiervar} together with Proposition \ref{decompprinciple}
is already enough to give a simple computation of the Dirac spectrum
of hyperbolic $n$-space $H^n$.
After removing a point $o$ from $H^n$ the space is isometric to a warped
product $S^{n-1} \times (0,\infty)$ where $S^{n-1}$ carries its standard
metric of constant sectional curvature 1 and $\rho(t) = \sinh(t)$.
By Lemma \ref{separiervar} the square of the Dirac operator on $H^n-\{ o \}$
is unitarily equivalent to $\bigoplus_{\mu\in spec(S^{n-1})} L_\mu$ where
$L_\mu = -\frac{d^2}{dt^2} + V_\mu(t)$, $V_\mu(t) = \frac{\mu \cosh(t) + \mu^2}
{\sinh(t)^2}$.
All Dirac eigenvalues $\mu$ of $S^{n-1}$ are nonzero.

Since $V_\mu(t) \to 0$ for $t \to \infty$ we know $spec_e(L_\mu) = [0,\infty)$
where $L_\mu$ is acting on the Hilbert space $L^2([1,\infty),\C,dt)$ with say
Dirichlet boundary conditions at $t=1$, see \cite[p.\ 1448, Thm.\ 16(b)]
{dunford-schwartz63a}.
By Proposition \ref{decompprinciple} we have that $spec_e((D^{H^n})^2) =
spec_e((D^{H^n-B_1(o)})^2) = [0,\infty)$.

One checks \cite{baier97a} that $D_\mu$ does not have any square integrable 
eigenfunctions on $(0,\infty)$, see also \cite{heinz82a,weidmann71a,
weidmann82a}.
In particular, there are no $L^2$-eigenspinors for the Dirac operator
on $H^n$, $spec_p(D^{H^n}) = \emptyset$.

We conclude $spec((D^{H^n})^2) = spec_c((D^{H^n})^2) = [0,\infty)$.
Finally, since $H^n$ is a simply connected symmetric space the spectrum of 
the Dirac operator is symmetric about 0.
In even dimensions this is automatic.
In odd dimensions the geodesic reflection about $o$ can be used to map
an eigenspinor or a Weyl sequence for $\lambda\in spec(D^{H^n})$ into one
for $-\lambda$.
We obtain
$$
spec(D^{H^n}) = spec_c(D^{H^n}) = \R.
$$
See \cite{bunke91a} for a computation of this spectrum using harmonic
analysis.
Note that there is an incorrect statement about the eigenvalue 0 in that
paper.
See also \cite{camporesi92a,camporesi-higuchi96a}.

\begin{definition}\punkt
Let $M$ be a hyperbolic manifold of finite volume.
Let $\E = N \times [0,\infty)$ be a cusp of $M$.
A spin structure of $M$ will be called {\bf trivial along} $\E$ if
the induced operator $\DN$ on $N$ has a nontrivial kernel, i.e.\
there exist $\vp \in \Cu(N,\SiM|_N)$, $\vp \not= 0$, but $\DN\vp = 0$.
\end{definition}

This terminology is justified by the fact that in the most prominent
case when $N$ is a flat torus,
the trivial (biinvariant) spin structure of $N$ is the only one
among its $2^{n-1}$ spin structures for which the Dirac operator has a 
nontrivial kernel.

The following theorem is our first main result.
It tells us that only two extremal cases can occur
for the type of spectrum of the Dirac operator on a hyperbolic
manifold of finite volume.
It can only be purely discrete spectrum or the whole real line.
It is the spin structure which is responsible for which of the
two cases occurs.

\begin{theorem}\punkt
\label{dichotomy}
Let $M$ be a hyperbolic manifold of finite volume equipped with
a spin structure.

If the spin structure is trivial along at least one cusp, then the Dirac 
spectrum satisfies
$$
spec(D) = spec_e(D) = \R .
$$

If the spin structure is nontrivial along all cusps, then the spectrum
is discrete,
$$
spec(D) = spec_d(D) .
$$
\end{theorem}

\noi
{\sc Proof.}
Recall decomposition (\ref{zerlegung}) of $M$ into a relatively
compact part and finitely many cusps,
$$
M = M_0 \udot \bigudot_{j=1}^k \E_j .
$$

We start with the case that the spin structure is trivial along at least
one end.
Hence $M$ has a cusp $\E_1 = N_1 \times [0,\infty)$ with metric $g_{\E_1} = 
e^{-2t}g_{N_1} + dt^2$ where $g_{N_1}$ is a flat metric and
$\ker(D^{N_1}) \not= 0$.
Choose $\vp \in \ker(D^{N_1})$, $\vp \not= 0$, such that 
$\nu\cdot\vp = i\cdot\vp$ or $\nu\cdot\vp = -i\cdot\vp$.
This is possible since $\nu^2 = -1$.
W.l.o.g.\ let $\nu\cdot\vp = i\cdot\vp$.

Let $\lambda \in \R$.
We look at spinors $\si$ on $\E_1$ of the form 
\begin{equation}
\si = \alpha\vp
\label{separat}
\end{equation}
where $\alpha : [0,\infty) \to \C$ is a smooth function.

For $\alpha_\lambda(t)=e^{((n-1)/2 - \lambda i)t}$, $\si_\lambda = 
\alpha_\lambda\vp$, we see using (\ref{diracgauss}) and $H=1$
\begin{eqnarray*}
D^{}\si_\lambda 
&=& 
\nu \left\{ D^{N_1} \sigma_\lambda - \frac{n-1}{2} H \sigma_\lambda + 
\nabla^{M}_{\nu} \sigma_\lambda \right\} \\
&=& 
\nu \left\{ 0 - \frac{n-1}{2} \alpha_\lambda\vp + \dot{\alpha}_\lambda\vp 
\right\} \\
&=&
-\lambda i \nu \cdot \si_\lambda \\
&=&
\lambda \si_\lambda .
\end{eqnarray*}
For $0 < a < b < \infty$ we denote $\E_{1,a,b} := N_1 \times
[a,b] \subset \E_1 \subset M$.
For a spinor $\si = \alpha\vp$ of the form (\ref{separat}) one easily 
computes the $L^2$-norm
$$
||\si||^2_{L^2(\E_{1,a,b})} = ||\vp||^2_{N_1} \cdot \int_a^b 
|\alpha(t)|^2 e^{-(n-1)t} dt .
$$
W.l.o.g.\ assume $||\vp||^2_{N_1} = 1$.
Then we obtain for $\si=\si_\lambda$ 
$$
||\si_\lambda||^2_{L^2(\E_{1,a,b})} = b-a .
$$
Now choose smooth functions $\psi_m : \R \to \R$ such that
\begin{itemize}
\item
$\psi_m \equiv 0$ on $\R - [m-2,2m+2]$,
\item
$\psi_m \equiv \frac{1}{\sqrt{m}}$ on $[m,2m]$,
\item
$0 \le \psi_m \le \frac{1}{\sqrt{m}}$ everywhere,
\item
$\left|\frac{d\psi_m}{dt}\right| \le \frac{1}{\sqrt{m}}$ everywhere.
\end{itemize}
We extend $\si_m := \psi_m \si_\lambda$ by zero to all of $M$ and compute
$$
||\si_m||^2_{L^2(M)} \ge
\frac{1}{m}||\si_\lambda||^2_{L^2(\E_{1,m,2m})} = 1 ,
$$
\begin{eqnarray*}
||(D^{} - \lambda)\si_m||^2_{L^2(M')} 
&=&
|| \nabla\psi_m \cdot \si_\lambda ||^2_{L^2(M)} \\
&=&
|| \nabla\psi_m \cdot \si_\lambda ||^2_{L^2(\E_{1,m-2,m}\cup 
\E_{1,2m,2m+2})} \\
&\le&
|| \nabla\psi_m ||^2_{L^\infty} \cdot ||\si_\lambda ||^2_{L^2(\E_{1,m-2,m}
\cup \E_{1,2m,2m+2})} \\
&\le&
\frac{4}{m} .
\end{eqnarray*}
Hence $(D^{} - \lambda)\si_m \to 0$.
For arbitrary square-integrable $\chi$ on $M$ we have
\begin{eqnarray*}
|(\chi,\si_m)_{L^2(M)}| 
&=&
|(\chi,\si_m)_{L^2(\E_{1,m-2,2m+2})}| \\
&\le&
||\chi||_{L^2(\E_{1,m-2,2m+2})} \cdot ||\psi_m||_{L^\infty} \cdot 
||\si_\lambda||_{L^2(\E_{1,m-2,2m+2})} \\
&\le&
\stackrel{\to 0}{\overbrace{||\chi||_{L^2(\E_{1,m-2,2m+2})}}} \cdot
\stackrel{\to 1}{\overbrace{\frac{1}{\sqrt{m}}\sqrt{m+4}}}.
\end{eqnarray*}
Thus
$$
\si_m \rightharpoonup 0.
$$

This shows $\lambda \in spec_e(D^{})$ and hence
$$
spec_e(D) = \R .
$$

Let us now put $M' := \bigudot_{j=1}^k \E_j$ and turn to the case that the 
spin structure is nontrivial along all cusps, i.e.\ $\ker(D^{\partial M'})= 0$.

By Lemma \ref{separiervar} the square of the Dirac operator on $M'$ is
unitarily equivalent to $\bigoplus_{\mu\in spec(D^{\del M'})}L_\mu$ where
$L_\mu = -\frac{d^2}{dt^2} + \mu e^t + \mu^2 e^{2t}$ is a 
Schr\"odinger operator on $L^2([0,\infty),\C,dt)$ with potential 
$V_\mu(t) = \mu e^t + \mu^2 e^{2t}$.
We impose Dirichlet boundary conditions.

Note that all $\mu\in spec(D^{\del M'})$ are nonzero.
Since $V_\mu \to \infty$ for $t \to \infty$ the classical theory of 
Weyl and Titchmarsh \cite[p.\ 1448, Thm.\ 16(a)]{dunford-schwartz63a} tells us
that the spectrum of $L_\mu$ is purely discrete, $spec_e(L_\mu) = \emptyset$.

To see $spec_e(\bigoplus_{\mu}L_\mu) = \emptyset$ we show
that only finitely many $\mu\in spec(D^{\del M'})$ contribute to the 
spectrum in a given compact interval $[-C,C] \subset \R$.

Only finitely many $\mu\in spec(D^{\del M'})$ satisfy $-1/2 < \mu <0$.
For all $\mu\in spec(D^{\del M'}) - (-1/2,0)$ we see
that $V_\mu(t) = \mu e^t + \mu^2 e^{2t} \ge \mu^2 - |\mu|$, hence
$spec(L_\mu) \subset [\mu^2 - |\mu|,\infty)$.
Since the $\mu\in spec(D^{\del M'})$ form a discrete set with $\mu^2
\to \infty$ there are only finitely many $\mu$ for which $spec(L_\mu)
\cap [-C,C] \not= \emptyset$.
Thus $spec(\bigoplus_{\mu}L_\mu)$ is discrete.

By Proposition \ref{decompprinciple} 
$$
spec_e(D^2) = spec_e((D^{M'})^2) = 
spec_e(\bigoplus_{\mu\in spec(D^{\del M'})}L_\mu) = \emptyset.
$$
Thus $spec_e(D) = \emptyset$ and the theorem is proven.
\qed

\begin{corollary}\punkt
\label{einende}
Let $M$ be a 2- or 3-dimensional hyperbolic manifold of finite volume 
equipped with a spin structure.
Let $M$ have exactly one cusp.

Then the spectrum of the Dirac operator is discrete,
$$
spec(D) = spec_d(D) .
$$
\end{corollary}

{\sc Proof.}
Decomposition (\ref{zerlegung}) of $M$ is in this case
$$
M = M_0 \udot \E .
$$
Here $M_0$ is a compact manifold with boundary $S^1$ or $T^2$ respectively.
It is well-known that the trivial spin structure on $S^1$ and $T^2$
do not bound a spin structure on a compact manifold.
Indeed, they generate spin cobordism $\Omega^{Spin}_1$ and $\Omega^{Spin}_2$
respectively, see e.g.\ \cite[p.\ 91]{lawson-michelsohn89a}.
Hence the spin structure must be nontrivial along $\E$.
Theorem \ref{dichotomy} yields the assertion.
\qed

\begin{corollary}\punkt
\label{immerdiskret}
Every 2- or 3-dimensional oriented hyperbolic manifold of finite volume
has a spin structure such that the spectrum of the Dirac operator 
is discrete,
$$
spec(D) = spec_d(D) .
$$
\end{corollary}

{\sc Proof.}
Again we look at decomposition (\ref{zerlegung}).
Chopping off the cusps yields the compact manifold $\barM_0$ with boundary.
The boundary is a disjoint union of $S^1$s or 2-tori.
Gluing in disks or solid tori we obtain an oriented closed manifold $M'$.
In dimension 2 and 3 all orientable manifolds are spin.
In three dimensions this follows from triviality of the tangent bundle.
Pick a spin structure on $M'$ and restrict it to $\barM_0$.
Since the trivial spin structures on $S^1$ and $T^2$ do not bound,
the induced spin structure must be nontrivial on all boundary
components.
Extending the spin structure to $M$ yields a spin structure which
is nontrivial along all cusps.
Hence Theorem \ref{dichotomy} yields the statement.
\qed

As we shall see in the last section a surface of finite volume with at
least two cusps can always be given a spin structure such that
$spec(D) = \R$.
Hence both cases in Theorem \ref{dichotomy} occur.
In three dimensions this depends on the manifold.
It can happen that the spectrum is always discrete even if the manifold
has more than one cusp.
If the hyperbolic manifold is given as the complement of a link in $S^3$,
then there is simple criterion to decide if there is a spin structure
such that $spec(D) = \R$.
This involves counting of overcrossings (Theorem \ref{verschlingung}).
See the last section for examples.


\section{Domain Monotonicity}

In order to study spectral degeneration in the next section we need
a tool known as domain monotonicity in the spectral theory of the
Laplace-Beltrami operator.
We have to find a version for the Dirac operator.
It will allow us to estimate the spectrum by decomposing the manifold
into pieces and controlling the spectrum of the individual pieces.
When doing this new boundary components appear and we have to exhibit
suitable boundary conditions.

Domain monotonicity can be conveniently expressed in terms of eigenvalue
counting functions.
Let $\barM$ be an $n$-dimensional compact Riemannian manifold with smooth 
boundary $\dM$.
Let $L$ be a formally self-adjoint elliptic differential operator acting
on sections of a Hermitian or Riemannian vector bundle defined over $\barM$.
Let the domain of $L$, specified by boundary conditions $\mathcal{B}$,
be such that $L$ becomes essentially self-adjoint.
Denote the corresponding self-adjoint extension by $\bar{L}$.
For any interval $I \subset \R$ we introduce the {\em eigenvalue counting 
function}
$$
\NN_{L,\barM}^\mathcal{B} (I) := \sharp (spec(\bar{L}) \cap I) .
$$
By passing to $n$-dimensional submanifolds of $\barM$ it is possible
to estimate $\NN_{L,\barM}^\mathcal{B} (I)$ from above and from below.
The two estimates are quite different in nature.
Let us start with the simpler one, the estimate from below.
Recall that an operator $d$ is called {\em overdetermined elliptic} if
its principal symbol $\si_d(\xi)$ is injective for all nonzero covectors
$\xi \in T^\ast M$.

\begin{proposition}\punkt
\label{monoton1}
{\em (Domain Monotonicity I)}
\newline
Let $\barM$ be a compact Riemannian manifold with smooth boundary
$\dM$.
Let $d:\Cu(M,E) \to \Cu(M,F)$ be an overdetermined elliptic linear
differential operator of first order, defined on Hermitian vector
bundles $E$ and $F$.
Put $L = d^\ast d$.
Let $\bar{N} \subset \barM$ be a compact submanifold with smooth boundary 
$\partial N$, $\dim(\bar{N})=\dim(\barM)$.
We impose Dirichlet boundary conditions, $\DD(L) = \{ \vp \in \Cu(\barM)\ |\
\vp|_{\dM} = 0 \}$, similarly for $\bar{N}$.

Then for any $x > 0$
$$
\NN_{L,\barM}^{\mbox{\em\tiny Dirichlet}} [0,x) \ge
\NN_{L,\bar{N}}^{\mbox{\em\tiny Dirichlet}} [0,x) .
$$
\end{proposition}

\begin{center}
\pspicture(1,0)(14,6)

\pscustom[linecolor=white,fillstyle=solid,fillcolor=lightgray]{
  \pscurve(8,1.3)(7.65,2)(7.5,3)(7.65,4)(8,4.7)
  \pscurve(8.01,4.7)(10,4.5)
  \pscurve(10.01,4.5)(10.01,1.5)
  \pscurve(10,1.5)(8,1.3)}

\psellipse(3.5,3)(0.5,1.5)
\psecurve[linecolor=white,linewidth=2pt](3,3)(3.5,4.5)(3.85,4)(4,3)(3.85,2)(3.5,1.5)(3,3)
\psecurve[linestyle=dotted](3,3)(3.5,4.5)(4,3)(3.5,1.5)(3,3)
\pscustom[fillstyle=solid,fillcolor=gray]{\psellipse(10,3)(0.5,1.5)}

\pscurve(3.5,1.5)(5,1.3)(6.5,1)(8,1.3)(10,1.5)
\pscurve(3.5,4.5)(5,4.7)(6.5,5)(8,4.7)(10,4.5)

\pscurve(5.8,3)(6.0,2.8)(6.4,2.7)(6.8,2.8)(7.0,3)
\pscurve(6.0,2.8)(6.4,3.1)(6.8,2.8)

\pscurve(8,1.3)(7.65,2)(7.5,3)(7.65,4)(8,4.7)

\rput(8.6,3){\psframebox*[framearc=0.5]{$\bar{N}$}}
\rput(4.4,5){$\bar{M}$}

\endpspicture
\end{center}
\begin{center}
\bf Fig.~5
\end{center}

{\sc Proof.}
The operator $\bar{L}$ is the self-adjoint operator associated with
the closed semi-bounded quadratic form 
$$
q(\vp) := (d\vp,d\vp)_{L^2}
$$
with a form core given by $\Cu_0(M,\SiM)$, cf.\
\cite[VIII.6]{reed-simon80a}.

Extension by zero yields an embedding $\Cu_0({N},\Si N) \hookrightarrow
\Cu_0(M,\SiM)$ and the quadratic form for $\bar{N}$ is simply the restriction
of the quadratic form for $\barM$.
The variational characterization of eigenvalues yields the proposition.
\qed

{\bf Example.}
If $d:\Cu(M,\R) \to \Cu(M,T^\ast M)$ is exterior differentiation, then
the proposition yields the standard domain monotonicity for the 
Laplace operator $L = \Delta$.

{\bf Example.}
If $d=D:\Cu(M,\SiM) \to \Cu(M,\SiM)$ is the Dirac operator on a Riemannian
spin manifold, then we obtain a monotonicity principle for the square
of the Dirac operator $L = D^2$.
We will use this in the next section.

For the reverse estimate we assume for simplicity that $M$ is a closed
manifold.
By a {\em decomposition} of $M$ we mean finitely many submanifolds
$\bar{M}_1 \bis \bar{M}_k$ of $M$ with smooth boundaries $\partial M_j$
and $\dim(M_j)=\dim(M)$ such that 
\begin{itemize}
\item
The interiors $\intM_j$ are pairwise disjoint,
\item
$M=\bigcup_{j=1}^k \barM_j$ .
\end{itemize}

\begin{center}
\pspicture(1,0)(14,6)

\pscustom[linecolor=white,fillstyle=solid,fillcolor=lightgray]{
  \pscurve(5,4.7)(3.5,4.5)(2.5,3)(3.5,1.5)(5,1.3)
  \pscurve(5.01,1.3)(4.89,1.5)(4.65,2)(4.5,3)(4.65,4)(5.01,4.7)}
  
\pscustom[linecolor=white,fillstyle=solid,fillcolor=lightgray]{
  \pscurve(8,1.3)(7.65,2)(7.5,3)(7.65,4)(8,4.7)
  \pscurve(8,4.7)(8.3,4.7)(8.8,4.7)(10,4.5)(11,3)(10,1.5)(8,1.3)}

\psccurve(3.5,1.5)(5,1.3)(6.5,1)(8,1.3)(10,1.5)(11,3)(10,4.5)(8,4.7)(6.5,5)(5,4.7)(3.5,4.5)(2.5,3)

\pscurve(5.8,3)(6.0,2.8)(6.4,2.7)(6.8,2.8)(7.0,3)
\pscurve(6.0,2.8)(6.4,3.1)(6.8,2.8)

\pscurve(8,1.3)(7.65,2)(7.5,3)(7.65,4)(8,4.7)
\pscurve(5,1.3)(4.65,2)(4.5,3)(4.65,4)(5,4.7)

\rput(3.6,3){\psframebox*[framearc=0.5]{$\bar{M}_1$}}
\rput(5.6,4){$\bar{M}_2$}
\rput(9.2,3){\psframebox*[framearc=0.5]{$\bar{M}_3$}}
\rput(4.4,5){${M}$}

\endpspicture
\end{center}
\begin{center}
\bf Fig.~6
\end{center}

If $\barM$ is a compact Riemannian manifold with smooth boundary 
$\dM$ and exterior unit normal field $\nu$ and $L = d^\ast d$ is as above, 
then we call the boundary condition
$$
(\si_{d^\ast}(\nu)d\vp)|_{\dM} = 0
$$
the {\em natural boundary conditions} for $L$.

\begin{proposition}\punkt
\label{monoton2}
{\em (Domain Monotonicity II)}
\newline
Let $M$ be a closed Riemannian manifold.
Let $d:\Cu(M,E) \to \Cu(M,F)$ be an overdetermined elliptic linear
differential operator of first order, defined on Hermitian vector
bundles $E$ and $F$.
Put $L = d^\ast d$.
Let $M=\bigcup_{j=1}^k \barM_j$ be a decomposition of $M$ as explained
above.

We impose natural boundary conditions for the $\barM_j$.
Let $L$ together with the natural boundary conditions form a regular elliptic
boundary value problem.

Then for any $x > 0$
$$
\NN_{L,M} [0,x) \le
\sum_{j=1}^k \NN_{L,\barM_j}^{\mbox{\em\tiny Natural}} [0,x) .
$$
\end{proposition}

{\sc Proof.}
If $\bar{N}$ is a compact Riemannian manifold and $L = d^\ast d$,
then $L$ with domain $\DD(L) = \{ \vp \in \Cu(\bar{N},\Si N)\ |\  
(\si_{d^\ast}(\nu)d\vp)|_{\partial N}=0 \}$ is essentially self-adjoint 
\cite[V.12]{taylor96a}.
Look at the closed semi-bounded quadratic form
$$
q(\vp) := (d\vp,d\vp)_{L^2}
$$
with form core $\Cu(\bar{N},\Si N)$.
The Green's formula
$$
(d^\ast d\vp,\psi)_{L^2(\bar{N},\Si\bar{N})} -
(d\vp,d\psi)_{L^2(\bar{N},\Si\bar{N})} =
\int_{\partial M} \langle \si_{d^\ast}(\nu) d\vp,\psi \rangle
$$
shows that the self-adjoint $A$ operator associated with $q$ has domain
\begin{eqnarray*}
\DD(A) &=&
\{ \vp \in \DD(q)\ |\ \exists \chi \in L^2(\bar{N},\Si\bar{N}):
(\psi,\chi)_{L^2} = (d\psi,d\vp)_{L^2} \\
&& \forall \psi \in \DD(q) \} \\
&\supset&
\{ \vp \in H^2(\bar{N},\Si\bar{N})\ |\ (\si_{d^\ast}(\nu)d\vp)|_{\partial N} 
= 0 \} \\
&=& \DD(\bar{L}) .
\end{eqnarray*}
Hence $\bar{L} \subset A$ and since both operators are self-adjoint
$\bar{L} = A$.
Therefore the eigenvalues of $L = d^\ast d$ can be computed
using the quadratic form $q$ with $\Cu(\bar{N},\Si\bar{N})$ as 
space of admissible test sections.

Returning to our closed manifold $M$ with the decomposition $M=\barM_1
\cup \ldots \cup \barM_k$ we look at the isometric embedding
$$
\Cu(M,\SiM) \subset L^2(M,\SiM) \to
\bigoplus_{j=1}^k \Cu(\bar{M}_j,\Si\bar{M}_j) \subset
\bigoplus_{j=1}^k L^2(\bar{M}_j,\Si\bar{M}_j) ,
$$
$$
\vp \mapsto (\vp|_{\barM_1} \bis \vp|_{\barM_k}) .
$$
Under this embedding the quadratic form $q$ corresponding to $L = 
d^\ast d$ on $M$ is the restriction of the orthogonal sum 
$q_1 \oplus \ldots \oplus q_k$ of the forms for $\barM_j$.
Now the variational characterization of eigenvalues completes the proof.
\qed

{\bf Example.}
If $d:\Cu(M,\R) \to \Cu(M,T^\ast M)$ is exterior differentiation, then
we obtain Neumann boundary conditions for the Laplace operator $L = \Delta$.
More generally, let $d = \nb : \Cu(M,E) \to \Cu(M,T^\ast M \otimes E)$
be a Riemannian connection.
Then the above monotonicity
 principle holds for the operator $L = \nb^\ast\nb$
with {\em Neumann boundary conditions}:
$$
0 = \si_{\nb^\ast}(\nu)\nb\vp = - \nb_\nu \vp .
$$

{\bf Example.}
If $d=D:\Cu(M,\SiM) \to \Cu(M,\SiM)$ is the Dirac operator on a Riemannian
spin manifold, then the monotonicity principle for the square
of the Dirac operator $L = D^2$ uses boundary conditions
$$
\nu\cdot D\vp|_{\dM} = 0
$$ 
or equivalently
$$
D\vp|_{\dM} = 0.
$$


\section{Degeneration in Two Dimensions}

Now we study the behavior of the spectrum under the degeneration
process described in the first section.
If the limit manifold has continuous spectrum we expect that
the eigenvalues of the compact manifolds accumulate in the 
degeneration.
We will see that this is true and compute the accumulation rate.
We start with the 2-dimensional case.

\begin{theorem}\punkt
\label{ausart2d}
Let $M_i$ be a sequence of closed hyperbolic surfaces converging
to a noncompact hyperbolic surface $M$ of finite volume.
Let each $M_i$ have exactly $k$ tubes with trivial spin structure around 
closed geodesics of length $\ell\uij$ tending to zero.
Hence $M$ has exactly $2k$ cusps along which the spin structure is trivial.
Let $x > 0$.

Then the eigenvalue counting function for the Dirac operator satisfies
for sufficiently small $\ell\uij$:
{\em $$
\NN_{D,{M_i}}(-x,x) = \frac{4x}{\pi}\sum_{j=1}^k \log(1/\ell\uij) 
+ \Oo_{x}(1).
$$}
\end{theorem}

Here $\Oo_{x}(1)$ denotes an error term bounded as a function of $i$
where the bound is allowed to depend on $x$.

{\sc Proof.}
To keep the notation simple we restrict ourselves to the case that
the $M_i$ have exactly one degenerating tube with either trivial
or nontrivial spin structure, hence $k=0$ or $k=1$.
Recall from Section 1 that the tube $T_i$ is isometric to
$S^1 \times [-R_i,R_i]$ with warped product metric $ds^2 = \ell_i^2
\cosh(t)^2 d\theta^2 + dt^2$, $\theta \in S^1=\R/\Z$, $t \in [-R_i,R_i]$,
$R_i = \log(1/\ell_i) + c_0$.

Choose a constant $c_1 = c_1(x)$ such that for all nonzero eigenvalues
$\mu$ of $D^{S^1}$ we have
\begin{equation}
e^{c_1-c_0} \cdot |\mu| \cdot (e^{c_1-c_0} \cdot |\mu| - 1) > x^2 .
\label{grossgenug}
\end{equation}
Put 
$$
\tilde{T}_i := S^1 \times [-R_i+c_1,R_i-c_1] \subset T_i \subset M_i
$$
and 
$$
\tilde{M}_{i,0} := M_{i,0} \cup \left(S^1 \times [-R_i,-R_i+c_1] \right)\cup
\left(S^1 \times [R_i-c_1,R_i]\right) \subset M_i.
$$
Then $M_i = \tilde{M}_{i,0} \cup \tilde{T}_i$.

The Lichnerowicz formula (\ref{lichnerowicz}) says in our case
$$
D^2 = \nb^\ast\nb + \frac{scal}{4} = \nb^\ast\nb - \eh.
$$
Proposition \ref{monoton2} yields
\begin{eqnarray*}
\NN_{D,{M_i}}(-x,x) &=&
\NN_{D^2,M_i}[0,x^2) \\
&=&
\NN_{\nb^\ast\nb - 1/2,M_i}[0,x^2) \\
&=&
\NN_{\nb^\ast\nb,M_i}[1/2,x^2+1/2) \\
&\le&
\NN_{\nb^\ast\nb,M_i}[0,x^2+1/2) \\
&\le&
\NN_{\nb^\ast\nb,\tilde{T}_i}^{\mbox{\tiny Neumann}}[0,x^2+1/2) +
\NN_{\nb^\ast\nb,\tilde{M}_{i,0}}^{\mbox{\tiny Neumann}}[0,x^2+1/2) \\
&=&
\NN_{D^2,\tilde{T}_i}^{\mbox{\tiny Neumann}}[-1/2,x^2) +
\NN_{D^2,\tilde{M}_{i,0}}^{\mbox{\tiny Neumann}}[-1/2,x^2) 
.
\end{eqnarray*}
All $\tilde{M}_{i,0}$ are diffeomorphic and the metrics converge in the
$\Cu$-topology to the metric of the limit surface.
Thus the eigenvalues also converge and therefore 
$\NN_{D^2,\tilde{M}_{i,0}}^{\mbox{\tiny Neumann}}[-1/2,x^2) = \Oo_x(1)$.

Using Lemma \ref{separiervar} we obtain
$$
\NN_{D^2,\tilde{T}_i}^{\mbox{\tiny Neumann}}[-1/2,x^2) =
\sum_{\mu\in spec(D^{S^1})} \NN_{L_\mu,[-R_i+c_1,R_i-c_1]}^{\mbox{\tiny Neumann}}[-1/2,x^2)
$$
where $L_\mu = -\frac{d^2}{dt^2} + V_\mu$, $V_\mu(t) = \mu \cdot
\frac{\ell_i\sinh(t) + \mu}{\ell_i^2\cosh(t)^2}$.
We can estimate the potential $V_\mu$ on $[-R+c_1,R-c_1] =
[-\log(1/\ell)-c_0+c_1,\log(1/\ell)+c_0-c_1]$ for nonzero $\mu \in
spec(D^{S^1})$ as follows:
\begin{eqnarray*}
V_\mu(t) &\ge& 
|\mu| \frac{|\mu| - \ell |\sinh(t)|}{\ell^2 \cosh(t)^2} \\
&\ge&
|\mu| \left(\frac{|\mu|}{\ell^2 \cosh(t)^2} - \frac{1}{\ell\cosh(t)}\right)\\
&=&
\frac{|\mu|}{\ell\cosh(t)}\left(  \frac{|\mu|}{\ell\cosh(t)} - 1 \right) \\
&\ge&
\frac{|\mu|}{\ell e^{R-c_1}}\left(  \frac{|\mu|}{\ell e^{R-c_1}} - 1 \right) \\
&=&
e^{c_1-c_0}|\mu|(e^{c_1-c_0}|\mu|-1) \\
&>&
x^2
\end{eqnarray*}
by (\ref{grossgenug}).
Hence for nonzero $\mu\in spec(D^{S^1})$ all eigenvalues of $L_\mu$
with Neumann boundary conditions are bigger than $x^2$, i.e.\ 
$\NN_{L_\mu,[-R_i+c_1,R_i-c_1]}^{\mbox{\tiny Neumann}}[-1/2,x^2) = 0$.
Denote the multiplicity of the eigenvalue $0$ in $spec(D^{S^1})$ by
$mult(0)$.
We have shown
\begin{eqnarray*}
\NN_{D,{M_i}}(-x,x) &\le&
\NN_{D^2,\tilde{T}_i}^{\mbox{\tiny Neumann}}[-1/2,x^2) +
\NN_{D^2,\tilde{M}_{i,0}}^{\mbox{\tiny Neumann}}[-1/2,x^2) \\
&=&
\NN_{D^2,\tilde{T}_i}^{\mbox{\tiny Neumann}}[-1/2,x^2) + \Oo_x(1) \\
&=&
mult(0)\cdot \NN_{-\frac{d^2}{dt^2},[-R_i+c_1,R_i-c_1]}^{\mbox{\tiny Neumann}}[-1/2,x^2) 
+ \Oo_x(1) \\
&=& 
mult(0)\cdot \NN_{-\frac{d^2}{dt^2},[-R_i+c_1,R_i-c_1]}^{\mbox{\tiny Neumann}}[0,x^2) 
+ \Oo_x(1) \\
&=& 
mult(0) \cdot \frac{x\cdot 2(R_i-c_1)}{\pi} + \Oo_x(1) \\
&=&
mult(0) \cdot \frac{2xR_i}{\pi} + \Oo_x(1) .
\end{eqnarray*}
In case the spin structure is nontrivial along $T_i$ we have $mult(0)=0$
and 
$$
\NN_{D,{M_i}}(-x,x) = \Oo_x(1).
$$
The theorem is proven in this case.
If the spin structure is trivial, we have $mult(0)=2$, hence
\begin{eqnarray*}
\NN_{D,{M_i}}(-x,x) &\le& \frac{4xR_i}{\pi} + \Oo_x(1) \\
&=& 
\frac{4x}{\pi}(\log(1/\ell_i)+c_0) + \Oo_x(1) \\
&=&
\frac{4x}{\pi}\log(1/\ell_i) + \Oo_x(1).
\end{eqnarray*}
In this case we also need a lower bound which is easily obtained by
applying Proposition \ref{monoton1} and Lemma \ref{separiervar}:
\begin{eqnarray*}
\NN_{D,{M_i}}(-x,x) &=&
\NN_{D^2,M_i}[0,x^2) \\
&\ge&
\NN_{D^2,T_i}^{\mbox{\tiny Dirichlet}}[0,x^2) \\
&=&
\sum_{\mu\in spec(D^{S^1})}\NN_{L_\mu,[-R_i,R_i]}^{\mbox{\tiny Dirichlet}}[0,x^2). \\
&\ge&
mult(0)\cdot\NN_{L_0,[-R_i,R_i]}^{\mbox{\tiny Dirichlet}}[0,x^2) \\
&=& 
2 \cdot \frac{x\cdot 2R_i}{\pi} + \Oo(1)\\
&=&
\frac{4x}{\pi} \log(1/\ell_i) + \Oo_x(1) .
\end{eqnarray*}
\qed


\section{Manifolds Foliated by Hypersurfaces}

In three dimensions we have the problem that the degenerating tube is not a 
warped product so that the simple separation of variables of Lemma 
\ref{separiervar} does not apply.
But the tube is foliated by flat 2-tori as described in Section 1.
In order to take advantage of this we derive a formula relating the
square of the Dirac operator on a manifold foliated by hypersurfaces
to normal derivatives and operators acting on the leaves.

In this paper we will only need Corollary \ref{tubenichttrivial} at the end 
of this section.
The reader may skip this section at a first reading and only come back to it
when needed.
We hope that the general formula in Proposition \ref{foliation} will also
be useful in other contexts.

Let $M$ be a Riemannian spin manifold of dimension $n$.
Let $M$ be foliated by oriented (hence spin) hypersurfaces $\{ N\}$.
Denote the unit normal field to the foliation by $\nu$, its shape
operator by $B$, $B(X)=-\nb_X\nu$, and its mean curvature function 
by $H:=\frac{1}{n-1}\Tr B$.

Let $\SiM$ be the spinor bundle on $M$.
Recall that $\SiM |_N$ is the spinor bundle of $N$ if $n$ is odd.
If $n$ is even, then $\SiM |_N$ coincides with the sum of two copies
of the spinor bundle of $N$.
Clifford multiplication with respect to $N$ is given by
$$
X \otimes \vp \to X \cdot \nu \cdot \vp
$$
where the dot ``$\cdot$'' denotes Clifford multiplication with respect to $M$.
Recall equation (\ref{diracgauss})
$$
-\nu\cdot \DM = \DN - \frac{n-1}{2}H + \naM_\nu .
$$

One can also relate $\naM$ to the spinorial Levi-Civita connection $\naN$ 
for $N$ by
\begin{equation}
\naM_X \vp = \naN_X \vp + \frac{1}{2}B(X)\cdot\nu\cdot\vp
\label{gauss}
\end{equation}
see e.g.\ \cite[Prop.~2.1]{baer96b}.

We need one more piece of notation.
Define
$$
\DB := \sum_{i=1}^{n-1} e_i \cdot\nu\cdot \naN_{B(e_i)}
= \sum_{i=1}^{n-1} B(e_i) \cdot\nu\cdot \naN_{e_i}
$$
If $B$ happens to be a multiple of the identity, $B = c\cdot$Id, then 
$$
\DB =  
c\sum_{i=1}^{n-1} e_i \cdot\nu\cdot \naN_{e_i} = c\DN .
$$

\begin{proposition}\punkt
\label{foliation}
Let $M$ be an $n$-dimensional Riemannian spin manifold with Ricci
curvature $Ric$.
Let $M$ be foliated by oriented hypersurfaces $\{ N \}$ as described above.
Then
\begin{eqnarray*}
(\DM)^2 &=& (\DN)^2 - (\naM_\nu)^2 + (n-1)H\naM_\nu
+ \naM_{\nb_\nu\nu}- \DB \\
&& - \frac{n-1}{2} (\naN H)\cdot\nu - \frac{(n-1)^2}{4}H^2
+\eh |B|^2 - \eh\nu \cdot Ric(\nu) .
\end{eqnarray*}
\end{proposition}

Here $|B|$ denotes the Hilbert-Schmidt norm of $B$, i.e.\ $|B|^2 = \sum_j
\lambda_j^2$ where $\lambda_1 \bis \lambda_{n-1}$ are the eigenvalues of $B$.

{\sc Proof.}
Let $e_1\bis e_{n-1}$ be a local
orthonormal tangent frame to one leaf of the foliation.
We locally solve the following linear ordinary differential equation
in the normal direction:
\begin{equation}
\nb_\nu e_j = - \< \nb_\nu \nu, e_j\>\nu .
\label{extend}
\end{equation}
Here $\nb$ denotes the Levi-Civita connection on $TM$.
We claim that this extends the frame to an orthonormal frame
$e_1\bis e_{n-1},\nu$ on an open subset of $M$.

Namely, we compute
\begin{eqnarray*}
\del_\nu \< \nu , e_j \> 
&=& \< \nb_\nu \nu , e_j \> + \< \nu , \nb_\nu e_j \> \\
&\stackrel{(\ref{extend})}{=}& 
\< \nb_\nu \nu , e_j \> + \< \nu , - \< \nb_\nu \nu, e_j\>\nu\> \\
&=& \< \nb_\nu \nu , e_j \> - \< \nb_\nu \nu, e_j\> \\
&=& 0
\end{eqnarray*}
and
\begin{eqnarray*}
\del_\nu \< e_i , e_j \>
&=& \< \nb_\nu e_i , e_j \> + \< e_i , \nb_\nu e_j \> \\
&\stackrel{(\ref{extend})}{=}&  
\< - \< \nb_\nu \nu, e_i\>\nu , e_j \> +\< e_i , - \< \nb_\nu\nu, e_j\>\nu\>\\
&=& 0.
\end{eqnarray*}

Let $\RSi$ be the curvature tensor on $\SiM$.
Recall \cite[Prop.~2.3]{baer96b} that Clifford multiplication 
by $\nu$ anticommutes with $\DN$
\begin{equation}
\DN(\nu\cdot\vp) = -\nu\cdot\DN\vp .
\label{anti}
\end{equation}

Now let us start the computation of $(\DM)^2$.
Squaring (\ref{diracgauss}) we obtain, using (\ref{anti}) and 
$\DN(f\vp) = (\naN f)\cdot \nu\cdot\vp + f\DN\vp$,
\begin{eqnarray}
(\DM)^2 &=& \left(\nu\cdot\DN - \frac{n-1}{2}H\nu + \nu\cdot\naM_\nu\right)
\left(\nu\cdot\DN - \frac{n-1}{2}H\nu + \nu\cdot\naM_\nu\right) \nonumber \\
&=& (\DN)^2 
- \frac{n-1}{2} \nu\cdot(\naN H)\cdot\nu\cdot\nu - \frac{n-1}{2}H\DN
+ \DN\naM_\nu \nonumber \\
&& + \frac{n-1}{2}H\DN - \frac{(n-1)^2}{4}H^2 + \frac{n-1}{2}H\naM_\nu
\nonumber \\
&& + \nu\cdot(\nb_\nu\nu)\cdot\DN - \naM_\nu\DN \nonumber \\
&& + \frac{n-1}{2}\cdot\del_\nu H - \frac{n-1}{2}H\nu\cdot(\nb_\nu\nu)
+ \frac{n-1}{2}H\naM_\nu  \nonumber \\
&& + \nu\cdot(\nb_\nu\nu)\cdot\naM_\nu - (\naM_\nu)^2 \nonumber \\
&=& (\DN)^2 - (\naM_\nu)^2
- \frac{n-1}{2} (\naN H)\cdot\nu
+ [\DN,\naM_\nu] \nonumber \\
&& + (n-1)H\naM_\nu + \nu\cdot(\nb_\nu\nu)\cdot\naM_\nu
+ \nu\cdot(\nb_\nu\nu)\cdot\DN \nonumber \\
&& - \frac{(n-1)^2}{4}H^2 + \frac{n-1}{2}\cdot\del_\nu H 
- \frac{n-1}{2}H\nu\cdot(\nb_\nu\nu)
\label{step1}
\end{eqnarray}
A simple computation yields the well-known formula
\begin{equation}
\sum_{j=1}^{n-1} e_j\cdot R^{\Si}(e_j,\nu) = \eh Ric(\nu) .
\label{ricci}
\end{equation}
We get
\begin{eqnarray}
[\DN,\naM_\nu] &=& \DN\naM_\nu - \naM_\nu\DN \nonumber \\
&=& \DN\naM_\nu - \naM_\nu\sum_{j=1}^{n-1}e_j\nu\naN_{e_j} \nonumber \\
&\stackrel{(\ref{extend})}{=}& 
\sum_{j=1}^{n-1}\<\nb_\nu\nu,e_j\>\nu\nu\naN_{e_j}
- \sum_{j=1}^{n-1}e_j(\nb_\nu\nu)\naN_{e_j}
+ \sum_{j=1}^{n-1}e_j\nu [\naN_{e_j},\naM_\nu ] \nonumber \\
&\stackrel{(\ref{gauss})}{=}& 
-\naN_{\nb_\nu\nu} + \sum_{j=1}^{n-1}(\nb_\nu\nu)e_j\naN_{e_j}
+2\sum_{j=1}^{n-1}\< e_j,\nb_\nu\nu\> \naN_{e_j} \nonumber \\
&&+ \sum_{j=1}^{n-1}e_j\nu [\naM_{e_j}-\frac{1}{2}B(e_j)\nu,\naM_\nu ] \nonumber \\
&=& (\nb_\nu\nu)\nu\DN + \naN_{\nb_\nu\nu}
+ \sum_{j=1}^{n-1}e_j\nu(\naM_{[e_j,\nu]}+\RSi(e_j,\nu)) \nonumber \\
&& + \frac{n-1}{2}H\naM_\nu
+ \frac{1}{2}\sum_{j=1}^{n-1}e_j\nu\naM_\nu\circ B(e_j)\nu \nonumber \\
&=& (\nb_\nu\nu)\nu\DN + \naN_{\nb_\nu\nu}\nonumber \\
&&+ \sum_{j=1}^{n-1}e_j\nu\naM_{-B(e_j)-\naM_\nu e_j}
- \nu\sum_{j=1}^{n-1}e_j\RSi(e_j,\nu) + \frac{n-1}{2}H\naM_\nu \nonumber \\
&&+ \naM_\nu\circ\frac{1}{2}\sum_{j=1}^{n-1}e_j\nu B(e_j)\nu
- \frac{1}{2}\sum_{j=1}^{n-1}(\naM_\nu e_j)\nu B(e_j)\nu \nonumber \\
&&- \frac{1}{2}\sum_{j=1}^{n-1}e_j(\nb_\nu\nu)B(e_j)\nu \nonumber \\
&\stackrel{(\ref{gauss})(\ref{extend})(\ref{ricci})}{=}& 
(\nb_\nu\nu)\nu\DN + \naN_{\nb_\nu\nu}\nonumber \\
&&-\DB +\sum_{j=1}^{n-1}e_j\nu\cdot\eh B(-B(e_j))\nu   
+ \sum_{j=1}^{n-1}e_j\nu\naM_{\<\nb_\nu\nu,e_j\>\nu}
- \eh\nu Ric(\nu) \nonumber \\
&&+ \frac{n-1}{2}H\naM_\nu - \naM_\nu\circ\frac{n-1}{2}H 
+ \frac{1}{2}\sum_{j=1}^{n-1}\<\nb_\nu\nu,e_j\>\nu\nu B(e_j)\nu\nonumber \\
&&+ \frac{1}{2}\sum_{j=1}^{n-1}(\nb_\nu\nu)e_jB(e_j)\nu
+ \sum_{j=1}^{n-1}\<\nb_\nu\nu,e_j\> B(e_j)\nu \nonumber \\
&=& (\nb_\nu\nu)\nu\DN + \naN_{\nb_\nu\nu} -\DB +\eh |B|^2
+ (\nb_\nu\nu)\nu\naM_{\nu}- \eh\nu Ric(\nu) \nonumber\\
&&+ \frac{n-1}{2}H\naM_\nu - \frac{n-1}{2}\del_\nu H - \frac{n-1}{2}H\naM_\nu
- \frac{1}{2} B(\nb_\nu\nu)\nu\nonumber \\
&&- \frac{n-1}{2}H(\nb_\nu\nu)\nu
+ B(\nb_\nu\nu)\nu \nonumber \\
&=&(\nb_\nu\nu)\nu\DN + \naN_{\nb_\nu\nu} - \DB +\eh |B|^2
+ (\nb_\nu\nu)\nu\naM_{\nu}- \eh\nu Ric(\nu) \nonumber \\
&& - \frac{n-1}{2}\del_\nu H + \frac{1}{2} B(\nb_\nu\nu)\nu
- \frac{n-1}{2}H(\nb_\nu\nu)\nu
\label{step2}
\end{eqnarray}
Plugging (\ref{step2}) into (\ref{step1}) yields
\begin{eqnarray}
(\DM)^2 &=& (\DN)^2 - (\naM_\nu)^2
- \frac{n-1}{2} (\naN H)\nu \nonumber \\
&& + (n-1)H\naM_\nu + \nu(\nb_\nu\nu)\naM_\nu
+ \nu(\nb_\nu\nu)\DN \nonumber \\
&& - \frac{(n-1)^2}{4}H^2 + \frac{n-1}{2}\del_\nu H 
- \frac{n-1}{2}H\nu(\nb_\nu\nu) \nonumber \\
&&+ (\nb_\nu\nu)\nu\DN + \naN_{\nb_\nu\nu} - \DB +\eh |B|^2
+ (\nb_\nu\nu)\nu\naM_{\nu}- \eh\nu Ric(\nu) \nonumber \\
&& - \frac{n-1}{2}\del_\nu H + \frac{1}{2} B(\nb_\nu\nu)\nu
- \frac{n-1}{2}H(\nb_\nu\nu)\nu \nonumber \\
&=&(\DN)^2 - (\naM_\nu)^2
- \frac{n-1}{2} (\naN H)\nu + (n-1)H\naM_\nu \nonumber \\
&& - \frac{(n-1)^2}{4}H^2 + \naN_{\nb_\nu\nu} - \DB +\eh |B|^2
- \eh\nu Ric(\nu) + \frac{1}{2} B(\nb_\nu\nu)\nu \nonumber \\
&\stackrel{(\ref{gauss})}{=}& 
(\DN)^2 - (\naM_\nu)^2 + (n-1)H\naM_\nu
+ \naM_{\nb_\nu\nu}- \DB \nonumber \\
&& - \frac{n-1}{2} (\naN H)\cdot\nu - \frac{(n-1)^2}{4}H^2
+\eh |B|^2 - \eh\nu \cdot Ric(\nu) 
\nonumber
\end{eqnarray}
\qed

Now let us specialize to the situation $M = N \times I$ where $N$ is a closed
$(n-1)$-dimensional spin manifold, $I \subset \R$ is an interval, and
$M$ carries a metric of the form
$$
ds^2 = g_r + dr^2
$$
where $g_r$ is a 1-parameter family of metrics on $N$.
The foliation is given by the leaves $N \times \{ r \}$.
Then $\nb_\nu\nu = 0$ and the formula in Proposition \ref{foliation}
simplifies to 
\begin{eqnarray*}
(\DM)^2 &=& (\DN)^2 - (\naM_\nu)^2 + (n-1)H\naM_\nu
- \DB \\
&& - \frac{n-1}{2} (\naN H)\cdot\nu - \frac{(n-1)^2}{4}H^2
+\eh |B|^2 - \eh\nu \cdot Ric(\nu) .
\end{eqnarray*}
We fix $r_0 \in I$.
Let $P_r$ denote parallel transport from $N \times \{ r_0 \}$ to
$N \times \{ r \}$ along $\nu$.
It is easy to see that 
$$
U : L^2(M,\SiM) \to L^2(I,L^2(N,g_{r_0},\SiM|_{N \times \{ r_0 \}}),dt),
$$
$$
(U\si)(r)(x) = \sqrt{\frac{dvol_{g_r}}{dvol_{g_{r_0}}}} P_r^{-1} \si(x,r),
$$ 
is a Hilbert space isometry.

\begin{corollary}\punkt
\label{foliationspezial}
The square of the Dirac operator on $M$, $(D^M)^2$, transforms under $U$
into the following Schr\"odinger operator acting on functions with values in 
$L^2(I,L^2(N,g_{r_0},\SiM|_{N \times \{ r_0 \}}),dt)$:
$$
U(D^M)^2U^{-1} = - \frac{d^2}{dr^2} + V
$$
where
$$
V = - \frac{n-1}{2}\frac{\del H}{\del r}
+ \frac{|B|^2}{2} + P_r^{-1}\left((D^N)^2 - \DB - \frac{n-1}{2} (\naN H)\nu 
- \eh\nu Ric(\nu)\right)P_r
$$
\end{corollary}

{\sc Proof.}
The first variation formula for the volume element of hypersurfaces tells
us
\begin{eqnarray}
\frac{d}{dr} \log \sqrt{\frac{dvol_{g_r}}{dvol_{g_{r_0}}}} &=&
\eh \frac{\frac{d}{dr}dvol_{g_r}}{dvol_{g_r}} \nonumber\\
&=&
- \eh \Tr(B) \nonumber\\
&=& 
- \frac{n-1}{2} H .
\label{varhyp}
\end{eqnarray}
We compute
\begin{eqnarray}
U\naM_\nu U^{-1} (v) &=&
U\naM_\nu \left( \sqrt{\frac{dvol_{g_{r_0}}}{dvol_{g_{r}}}} P_r v \right)
\nonumber\\
&=& 
U \left( \frac{n-1}{2} H \sqrt{\frac{dvol_{g_{r_0}}}{dvol_{g_{r}}}} P_r v
+ \sqrt{\frac{dvol_{g_{r_0}}}{dvol_{g_{r}}}} P_r \dot{v} \right)
\nonumber\\
&=&
\frac{n-1}{2} H v + \dot{v}
\label{normalAbl1}
\end{eqnarray}
and
\begin{eqnarray}
U(\naM_\nu)^2 U^{-1} (v) 
&=&
U\naM_\nu U^{-1} \left(\frac{n-1}{2} H v + \dot{v}\right) \nonumber\\
&=&
\frac{n-1}{2} H \left(\frac{n-1}{2} H v + \dot{v}\right) +
\frac{n-1}{2} \frac{\del H}{\del r} v + \frac{n-1}{2} H \dot{v}
+ \ddot{v} \nonumber\\
&=&
\ddot{v} + (n-1) H \dot{v} +
\left(\frac{(n-1)^2}{4} H^2 + \frac{n-1}{2} \frac{\del H}{\del r} \right) v .
\label{normalAbl2}
\end{eqnarray}
Equations (\ref{normalAbl1}) and (\ref{normalAbl2}) yield
$$
U\left( -(\naM_\nu)^2 + (n-1)H\naM_\nu - \frac{(n-1)^2}{4}H^2\right)U^{-1}(v)
$$
\begin{eqnarray}
&=&
-\ddot{v} - (n-1) H \dot{v} -
\left(\frac{(n-1)^2}{4} H^2 + \frac{n-1}{2} \frac{\del H}{\del r} \right) v
\nonumber\\
&& + \frac{(n-1)^2}{2} H^2 v + (n-1)H\dot{v} - \frac{(n-1)^2}{4}H^2
\nonumber\\
&=&
-\ddot{v} - \frac{n-1}{2} \frac{\del H}{\del r} v .
\end{eqnarray}
The corollary now follows from Proposition \ref{foliation}.
\qed

{\bf Example.}
Let us now look at the example of main interest in this paper, the 
tube around a closed geodesic in a hyperbolic 3-manifold.
Recall from Section 1 that $M = T[1,R]$ is isometric to $T^2 \times
[1,R]$ with Riemannian metric $ds^2 = g_r + dr^2$ where $g_r$ is
the flat metric on $T^2$ given by the lattice $\Gamma_r \subset \R^2$ 
spanned by the vectors $(2\pi\sinh(r),0)$ and $(\alpha_{i,j}\sinh(r), 
\ell_{i,j}cosh(r))$.

The shape operator $B$ has the eigenvalues $\tanh(r)$ and $\coth(r)$.
Hence $|B|^2 = \tanh(r)^2 + \coth(r)^2$, the mean curvature 
$H = \eh (\tanh(r) + \coth(r))$ is constant along the leaves and
$\nb^N H = 0$.
Since the sectional curvature of $T[r_1,r_2]$ is $-1$ we have
$Ric = -2\cdot\Id$.
Therefore Corollary \ref{foliationspezial} gives
$$
U(D^M)^2U^{-1} = -\frac{d^2}{dr^2} + \tanh(r)^2 + \coth(r)^2 - 2
+ P_r^{-1} \left( (D^N)^2 - \DB \right) P_r .
$$
It will be important to estimate the potential 
$$
V =\tanh(r)^2 + \coth(r)^2 - 2 + P_r^{-1} \left( (D^N)^2 - \DB \right) P_r 
$$
of this Schr\"odinger operator from below.
Note that $\DB$ is formally self-adjoint because $B$ is parallel along
the leaves.
If $\vp$ is a spinor field along a leaf, then
\begin{eqnarray*}
|\DB\vp| &=& \left|\sum_{j=1}^2 B(e_j) \cdot\nu\cdot \naN_{e_j}\vp\right| \\
&\le&
2\cdot |B|\cdot |\naN\vp|,
\end{eqnarray*}
hence
\begin{eqnarray*}
||\DB\vp||_{L^2(N,\SiM|_N)}^2 &\le& 
4\cdot |B|^2\cdot  ||\naN\vp||_{L^2(N,\SiM|_N)}^2\\
&=&
4\cdot |B|^2\cdot  ((\naN)^\ast\naN\vp,\vp)_{L^2(N,\SiM|_N)} \\
&=&
4\cdot |B|^2\cdot \left(\left((D^N)^2 - \frac{scal_N}{4}\right)\vp,
\vp\right)_{L^2(N,\SiM|_N)} \\
&=&
4 (\tanh(r)^2 + \coth(r)^2) \cdot  ||D^N\vp||_{L^2(N,\SiM|_N)}^2 \\
&\le&
4 (1 + \coth(1)^2) \cdot  ||D^N\vp||_{L^2(N,\SiM|_N)}^2 \\
&\le&
16 \cdot  ||D^N\vp||_{L^2(N,\SiM|_N)}^2 .
\end{eqnarray*}
Thus if $||D^N\vp||_{L^2(N,\SiM|_N)}^2 \ge \mu^2\cdot ||\vp||_{L^2(N,
\SiM|_N)}^2$, $\mu \ge 4$, we get
\begin{eqnarray}
\left(((D^N)^2 - \DB)\vp,\vp\right)_{L^2(N,\SiM|_N)} &\ge&
||D^N\vp||_{L^2(N,\SiM|_N)}^2 \nonumber\\
&& - 4\cdot ||D^N\vp||_{L^2(N,\SiM|_N)} \nonumber\\
&\ge&
\mu(\mu - 4) ||\vp||_{L^2(N,\SiM|_N)}^2 . \nonumber
\end{eqnarray}
Hence on each Hilbert subspace of $L^2(N,g_{r_0},\SiM|_{N \times \{ r_0 \}})$
which is left invariant by $P_r^{-1}(D^N)^2 P_r$ and by $P_r^{-1}\DB P_r$ 
on which $P_r^{-1}(D^N)^2 P_r \ge \mu^2$, $\mu \ge 4$, we know that
$P_r^{-1} ((D^N)^2-\DB) P_r \ge \mu(\mu - 4)$ and hence
\begin{equation}
V \ge \tanh(r)^2 + \coth(r)^2 - 2 + \mu(\mu - 4) > \mu(\mu - 4) - 1
\label{potentialvonunten}
\end{equation}
In order to proceed we need to control the eigenvalues of $(D^N)^2$.

\begin{lemma}\punkt
The smallest nonzero eigenvalue of $(D^N)^2$ on $N \times \{ r \}$
is monotonically decreasing in $r \in [1,R]$.
\end{lemma}

{\sc Proof.}
The Dirac eigenvalues of a flat torus $T^2 = \R^2 / \Gamma$ can be 
computed in terms of the dual lattice $\Gamma^\ast$.
A spin structure corresponds to a pair $\delta = (\delta_1,\delta_2)$,
$\delta_j = 0,1$, and the square of the Dirac operator for the corresponding 
spin structure has the eigenvalues
$$
4\pi^2\left| v - \eh(\delta_1 v_1 + \delta_2 v_2)\right|^2 
$$
where $v_1,v_2$ are a basis of $\Gamma^\ast$ and $v$ runs through 
$\Gamma^\ast$, cf.\ \cite{friedrich84a}.
In our case $\Gamma$ has the basis
$$
w_1 = \left(
\begin{array}{c}
2\pi\sinh(r) \\
0
\end{array}
\right)
,
w_2 = \left(
\begin{array}{c}
\alpha\sinh(r) \\
\ell\cosh(r)
\end{array}
\right) ,
$$
compare Section 1.
Hence a basis for $\Gamma^\ast$ is given by
$$
v_1 = \left(
\begin{array}{c}
\frac{1}{2\pi\sinh(r)} \\
-\frac{\alpha}{2\pi\ell\cosh(r)}
\end{array}
\right)
,
v_2 = \left(
\begin{array}{c}
0 \\
\frac{1}{\ell\cosh(r)}
\end{array}
\right) .
$$
Thus the eigenvalues
$$
\frac{(k_1 - \delta_1/2)^2}{\sinh(r)^2} +
\frac{(2\pi (k_2 - \delta_2/2) - \alpha (k_1 - \delta_1/2))^2}
{\ell^2\cosh(r)^2} ,
$$
$k_1, k_2 \in \Z$, are monotonically decreasing functions.
\qed

The lemma together with (\ref{potentialvonunten}) immediately implies

\begin{corollary}\punkt
\label{tubenichttrivial}
If $M = T[1,R]$ carries a nontrivial spin structure and if the smallest 
eigenvalue $\mu^2$ of $(D^N)^2$ on $N \times \{ R \}$ satisfies $\mu \ge 4$, 
then $(D^M)^2$ is unitarily equivalent to a Schr\"odinger operator
$$
-\frac{d^2}{dr^2} + V ,
$$
acting on Hilbert space-valued functions on $[1,R]$ with
$$
V \ge \mu(\mu - 4) - 1 .
$$
\end{corollary}


\section{Degeneration in Three Dimensions}

With the preparations in the previous section we are able to modify
the proof of Theorem \ref{ausart2d} such that it also works in three
dimensions.
In contrast to the 2-dimensional case there is no accumulation of
eigenvalues.

\begin{theorem}\punkt
\label{ausart3d}
Let $M_i$ be a sequence of closed hyperbolic 3-manifolds converging
to a noncompact hyperbolic 3-manifold $M$ of finite volume.
Let each $M_i$ have exactly $k$ tubes around 
closed geodesics of length $\ell\uij$ tending to zero.
Hence $M$ has exactly $k$ cusps.
Let $x > 0$.

Then the spin structure is nontrivial along all tubes and
the eigenvalue counting function for the Dirac operator remains bounded:
{\em $$
\NN_{D,{M_i}}(-x,x) = \Oo_{x}(1).
$$}
\end{theorem}

{\sc Proof.}
Again we restrict ourselves to the case that there is exactly one degenerating
tube, i.e.\ $k=1$.
Recall the decomposition of the manifolds
$$
M_i = M_{i,0} \udot T_i[0,R_i],
$$
$R_i = \eh\log(1/\ell_i)+c_0$.
Since $\dM_{i,0} = \del T_i[0,R_i]$ bounds the solid 2-torus $T_i[0,R_i]$
the induced spin structure on $\del T_i[0,R_i]$ must be nontrivial.

Look at decomposition (\ref{zerlegung}) of the limit manifold
$$
M = M_0 \udot \E
$$
where the cusp $\E = N \times [0,\infty)$ carries the warped product 
metric $e^{-2t}\cdot g_N + dt^2$.
Since we assumed compatibility of the spin structures of the $M_i$ and of
$M$, the spin structure of $M$ must also be nontrivial along $\E$.

Let $\mu_0$ be the smallest positive eigenvalue of the Dirac operator
on $(N,g_N)$.
Choose a constant $c_1 = c_1(x)$ such that
\begin{equation}
e^{c_1}\mu_0(e^{c_1}\mu_0 -4) -1 > x^2
\label{gutewahl}
\end{equation}
and
$$
e^{c_1}\mu_0 > 4.
$$
The number $e^{c_1}\cdot\mu_0$ is the smallest positive eigenvalue of the
Dirac operator on $(N, e^{-2c_1}g_N)$.
Put
$$
\tilde{M}_{i,0} := M_{i,0} \cup T_i[R_i - c_1,R_i] .
$$
This yields the following decomposition of the manifolds $M_i$:
$$
M_i = \tilde{M}_{i,0} \cup T_i[1,R_i - c_1] \cup T_i[0,1].
$$
Using Proposition \ref{monoton2} and the Lichnerowicz formula 
(\ref{lichnerowicz}) we get as in the proof of Theorem \ref{ausart2d}
\begin{eqnarray}
\NN_{D,{M_i}}(-x,x) &\le& 
\NN_{D^2,\tilde{M}_{i,0}}^{\mbox{\tiny Neumann}}[-3/2,x^2) +
\NN_{D^2,{T}_i[1,R_i - c_1]}^{\mbox{\tiny Neumann}}[-3/2,x^2) \nonumber\\
&& + \NN_{D^2,{T}_i[0,1]}^{\mbox{\tiny Neumann}}[-3/2,x^2) .
\label{ober1}
\end{eqnarray}
The same argument as in the proof of Theorem \ref{ausart2d} gives
\begin{equation}
\NN_{D^2,\tilde{M}_{i,0}}^{\mbox{\tiny Neumann}}[-3/2,x^2) =
\Oo_x(1) .
\label{ober2}
\end{equation}
The universal covering of the 1-tube $T_i[0,1]$ around the closed geodesic
$\gamma_i$ is the 1-tube $\mathcal{T}$ around a geodesic $\gamma$ in 
hyperbolic 3-space.
The group of deck transformations is isomorphic to $\Z$ and is generated
by a shift of length $\ell = \ell_i$ along $\gamma$ while rotating about the
angle $\alpha=\alpha_i$.
Denote this isometry of $\mathcal{T}$ by $A_{\ell,\alpha}$.
Note that $\alpha$ takes values in the compact interval $[-\pi,\pi]$.
As long as $\ell$ also takes values in a compact interval, say $\ell \in
[1/2,1]$, all eigenvalues vary in a bounded range and 
$\NN_{D^2,\mathcal{T}/A_{\ell,\alpha}}^{\mbox{\tiny Neumann}}[-3/2,x^2) 
= \Oo_x(1)$.

Now if $0 < \ell < \eh$ choose $m \in \N$ such that $m\ell \in [1/2,1]$.
Then $\mathcal{T}/A_{\ell,\alpha}$ is covered by $\mathcal{T}/
A_{m\ell,\alpha'}$.
Hence every eigenvalue of $\mathcal{T}/A_{\ell,\alpha}$ is also an eigenvalue
of $\mathcal{T}/A_{m\ell,\alpha'}$.
Therefore $\NN_{D^2,\mathcal{T}/A_{\ell,\alpha}}^{\mbox{
\tiny Neumann}}[-3/2,x^2) \le \NN_{D^2,\mathcal{T}/A_{m\ell,\alpha'}}^{\mbox{
\tiny Neumann}}[-3/2,x^2) = \Oo_x(1)$.
This shows $\NN_{D^2,\mathcal{T}/A_{\ell,\alpha}}^{\mbox{\tiny Neumann}}
[-3/2,x^2) = \Oo_x(1)$ for all $\alpha\in [-\pi,\pi]$ and $\ell\in (0,1]$.
Hence
\begin{equation}
\NN_{D^2,{T}_i[0,1]}^{\mbox{\tiny Neumann}}[-3/2,x^2) = \Oo_x(1).
\label{ober3}
\end{equation}
It remains to estimate $\NN_{D^2,{T}_i[1,R_i - c_1]}^{\mbox{
\tiny Neumann}}[-3/2,x^2)$.
By Corollary \ref{tubenichttrivial} the operator $D^2$ on ${T}_i[1,R_i - c_1]$ 
is unitarily equivalent to a Schr\"odinger operator $-\frac{d^2}{dr^2}+V_i$.
For sufficiently large $i$ the potential $V_i$ is bounded from below by 
$e^{c_1}\mu_0(e^{c_1}\mu_0 -4) -1$.
This follows from (\ref{gutewahl}) because the eigenvalues of $\del 
\tilde{M}_{i,0}$ converge to those of $\del\tilde{M}_{i}$.
We conclude
\begin{equation}
\NN_{D^2,{T}_i[1,R_i - c_1]}^{\mbox{\tiny Neumann}}[-3/2,x^2) = 0
\label{ober4}
\end{equation}
for sufficiently large $i$.
Plugging (\ref{ober2}), (\ref{ober3}), and (\ref{ober4}) into (\ref{ober1})
we obtain
$$
\NN_{D,{M_i}}(-x,x) = \Oo_x(1).
$$
\qed


\section{Spin Structures on Hyperbolic Manifolds}

The previous discussion has shown that the spectrum of the Dirac
operator depends in a crucial way on the spin structure.
This is true for the degeneration as well as for the $L^2$-spectrum
of a hyperbolic manifold of finite volume.
The fact that there is no spectral accumulation in three dimensions
has a topological reason.
Tubes necessarily carry a nontrivial spin structure because the trivial one
on the 2-torus does not bound.
In this last section we will discuss the question which kind of spin 
structures are actually carried by 2- or 3-dimensional hyperbolic manifolds 
of finite volume.

All (co-)homology groups in this section are to be taken with coefficients
$\zz$.
Recall that $H^1(M)$ acts simply transitively on the set of spin structures
of a spin manifold $M$.

{\bf The 2-dimensional case.}
Let $M$ be an oriented surface with $k$ ends.
Topologically $M$ is a closed surface $\bar{M}$ with $k$ points $p_1 \bis p_k$
removed.
Let $D_j$ denote small disks around $p_j$.
The Mayer-Vietoris sequence for the pair $(M, \bigudot_{j=1}^k D_j)$ yields
an exact sequence
\begin{equation}
0 \longrightarrow H^1(\bar{M}) \longrightarrow H^1(M) \longrightarrow
\bigoplus_{j=1}^k H^1(D_j - \{ p_j \}) \longrightarrow H^2(\bar{M}) 
\longrightarrow 0.
\label{MV1}
\end{equation}
Pick a spin structure on $\barM$ to identify spin structures with elements
of $H^1(\barM)$.
Take the restriction of this spin structure to $M$ and identify the spin
structures on $M$ with elements of $H^1(M)$.
The unique spin structure on $D_j$ induces the nontrivial spin structure
on $D_j - \{ p_j \} \simeq S^1$.
Hence (\ref{MV1}) tells us that the restriction mapping from spin structures
on $\barM$ to $M$ is injective and a spin structure on $M$ extends to $\barM$
if and only if it is nontrivial along all ends.

If we identify $H^1(D_j - \{ p_j \}) \cong \zz$ and 
$H^2(\barM) \cong \zz$, then the map $\bigoplus_{j=1}^k H^1(D_j - \{ p_j \}) 
\to H^2(\bar{M})$ corresponds to 
$(\zz)^k \stackrel{(1 \bis 1)}{\longrightarrow} \zz$.
Hence any spin structure on $M$ must be trivial along an even number of
ends.

To summarize, {\em a spin structure on $M$ corresponds uniquely to a spin 
structure on $\barM$ together with a choice of an even number of ends along
which the spin structure is trivial.}
$$
\left\{ \begin{array}{c}
        \mbox{spin structures}\\
        \mbox{on $M$} 
        \end{array} \right\}
\stackrel{1:1}{\longleftrightarrow}
\left\{ \begin{array}{c}
        \mbox{spin structures}\\
        \mbox{on $\barM$} 
        \end{array} \right\}
\times
\left\{ \begin{array}{c}
        \mbox{choices of an even}\\
        \mbox{number of ends of $M$}
        \end{array} \right\}
$$

In particular, on hyperbolic surfaces of finite area with more than one cusp
both cases in Theorem \ref{dichotomy} do occur.

To discuss Theorem \ref{ausart2d} let now $M$ be a closed oriented
surface.
Let $T \approx S^1 \times I$ be a tube around around a closed geodesic
$\gamma$.
The tube $T$ may carry the trivial or the nontrivial spin structure.
Is it possible to ``flip'' the spin structure, i.e.\ are both spin structures
on $T$ induced by some spin structure on $M$?

One can flip the spin structure if and only if there exists a cohomology class
in $H^1(M)$ acting nontrivially on $[\gamma] \in H_1(M)$, i.e.\ if and only if
the homology class $[\gamma]$ is nonzero in $H_1(M)$.
This is the case if and only if removing $\gamma$ does not decompose the 
surface into two connected components.

If $M-\gamma$ is disconnected which spin structure does $T$ carry?
Both connected components can be given a spin structure which is
nontrivial along all ends, cf.\ the discussion above.
These spin structures can be glued together to give a spin structure on
the original surface $M$.
Hence $T$ must carry the nontrivial spin structure in this case.
We note:

{\em The tube can carry both spin structures if and only if cutting along
$\gamma$ does not decompose the surface into two connected components.
In this case spectral accumulation in Theorem \ref{ausart2d} may or
may not occur depending on the choice of spin structure.
If $M-\gamma$ disconnects, then the tube carries the nontrivial spin structure
and does not contribute to the spectral accumulation.}

{\bf The 3-dimensional case.}
The proof of Theorem \ref{ausart3d} has shown that all tubes in a closed
hyperbolic 3-manifold carry the nontrivial spin structure.
This is responsible for the fact that there is no spectral accumulation
in three dimensions.
Spin structures on hyperbolic 3-manifolds of finite volume which are
trivial along some cusps do not occur as limits of spin structures on
closed hyperbolic 3-manifolds.
Do they exist at all?

Let us first show that just like in two dimensions {\em any spin structure on
a hyperbolic 3-manifold of finite volume is trivial along an even number
of cusps.}

Let $M = M_0 \udot \bigudot_{j=1}^k \E_j$ be a hyperbolic spin 3-manifold with
$k$ cusps.
Let the spin structure be trivial along $k_1$ cusps and nontrivial along
$k_2$ cusps, $k = k_1 + k_2$.
Chop off the ends to obtain the compact manifold $\barM_0$ with boundary.
Two of the three nontrivial spin structures on the 2-torus bound spin 
structures on the solid torus $S$, $T^2 = \del S$.
The third one can be transformed by some automorphism of $T^2$ into
one which bounds a spin structure on the solid torus.
Hence using appropriate gluing maps we can glue in solid tori to the boundary 
components of $\barM_0$ on which the spin structure is nontrivial and extend 
the spin structure.
We obtain a compact spin manifold $M_0'$ whose boundary consists of
$k_1$ tori.
The induced spin structure is trivial on all these boundary components.

Assume $k_1$ were odd, $k_1 = 2m+1$.
Choose $m$ pairs of boundary tori and identify them.
Since the spin structures on the tori are all trivial they can be glued
together.
We obtain a compact spin manifold $M_0''$ whose boundary consists of the
one remaining 2-torus.
The induced spin structure on this torus is trivial.
This contradicts the fact that the trivial spin structure on $T^2$ does
not bound.
\qed

Here is a criterion for when a boundary torus can inherit the trivial
spin structure.

\begin{lemma}\punkt
\label{randtorus}
Let $M$ be an oriented 3-manifold with boundary.
Let $T$ be a connected component of the boundary diffeomorphic to a
2-torus.
Then the following two assertions are equivalent:
\begin{itemize}
\item
$M$ carries a spin structure inducing the trivial spin structure on $T$.
\item
The inclusion map $T \hookrightarrow M$ induces an {\bf injective} map on 
the first homology
$$
H_1(T) \to H_1(M).
$$
\end{itemize}
\end{lemma}

{\sc Proof.}
A solid torus $S$ induces exactly two of the three nontrivial spin structures
on its boundary $\del S = T^2$.
Denote the three nontrivial spin structures on $T^2$ by $\SSS_1, \SSS_2, \SSS_3$.
Choose a spin structure on $M \cup_{(\Id,T^2)} S$.
Since the induced spin structure on $T^2$ bounds a spin structure on $S$
it is nontrivial, say $\SSS_1$.

The automorphisms of $T^2$ act transitively on $\{\SSS_1,\SSS_2,\SSS_3\}$.
Choose an automorphism $\Phi$ of $T^2$ such that $\Phi^\ast\SSS_1$
is the nontrivial spin structure on $T^2$ which is {\bf not} induced
by one on $S$.
Pick a spin structure on $M \cup_{(\Phi,T^2)} V$.
The induced spin structure on $T^2$ is again nontrivial but $\not=\SSS_1$.

We have found two spin structures on $M$ inducing two different nontrivial
spin structures on $T$.

{\bf Case 1.}
$H_1(T) \to H_1(M)$ is injective, i.e.\ $H^1(M) \to H^1(T)$ is surjective.

In this case $H^1(M)$ acts transitively on the spin structures of $T$
and in particular the trivial spin structure occurs.

{\bf Case 2.}
$H_1(T) \to H_1(M)$ is not injective, i.e.\ $H^1(M) \to H^1(T)$ is not
surjective.

In this case $\dim Im H^1(M) \le 1$ and hence $\sharp Im H^1(M) \le 2$
where $Im H^1(M)$ denotes the image of $H^1(M)$ in $H^1(T)$.
Therefore at most two spin structures are induced on $T$.
But as we have seen above there are two nontrivial spin structures
which do occur.
Hence $\dim Im H^1(M) = 1$ and $T$ inherits exactly two spin structures
both nontrivial.
\qed

A main source of hyperbolic manifolds of finite volume is given by complements
of links in $S^3$.
For such manifolds Lemma \ref{randtorus} can be translated in a very simple
criterion.

\begin{theorem}\punkt
\label{verschlingung}
Let $K \subset S^3$ be a link, let $M = S^3 - K$ carry a hyperbolic metric
of finite volume.

If the linking number of all pairs of components $(K_i,K_j)$ of $K$
is even,
$$
Lk(K_i,K_j) \equiv 0 \mbox{ mod } 2,
$$
$i \not= j$, then the spectrum of the Dirac operator on $M$ is discrete for all
spin structures,
$$
spec(D) = spec_d(D).
$$

If there exist two components $K_i$ and $K_j$ of $K$, $i \not= j$, with
odd linking number, then $M$ has a spin structure such that the spectrum of the
Dirac operator satisfies
$$
spec(D) = \R.
$$
\end{theorem}

{\sc Proof.}
Each component $K_j$ of $K$ corresponds to one cusp of $M$.
Let $K_1 \bis K_k$ be the components of $K$ and let $S_1 \bis S_k$ denote
thin solid tori around the link components.
The solid tori have to be pairwise disjoint.
Denote the boundary tori by $T_j = \del S_j$.

The Mayer-Vietoris sequence for the pair $(M, \bigudot_{j=1}^k S_j)$
yields an exact sequence
\begin{equation}
0 \longrightarrow \bigoplus_{j=1}^k H_1(T_j) \longrightarrow
H_1(M) \oplus \bigoplus_{j=1}^k H_1(S_j) \longrightarrow 0
\label{MV2}
\end{equation}
Choose a basis $\alpha_j,\beta_j$ of $H_1(T_j)$ such that
$\alpha_j$ generates the kernel of $H_1(T_j) \to H_1(S_j)$ and 
$\beta_j$ is represented by curve unlinked to the soul of $S_j$.
From (\ref{MV2}) we see that the map $\bigoplus_{j=1}^k H_1(T_j) \to
H_1(M)$ restricted to the span of $\alpha_1 \bis \alpha_k$ is injective.

Now let $c_j$ denote the linking numbers of $K_1$ and $K_j$, $j \ge 2$.
Then $\beta_1$ is homologous to $\sum_{j=2}^k c_j \alpha_j \in H_1(M)$,
see Figure 7.

\begin{center}
\pspicture(1,0)(14,7)

\psellipse[linewidth=2pt](7,3.5)(4,2)
\pscustom[linewidth=1pt,fillstyle=solid,fillcolor=lightgray]{
  \psellipse(7,3.5)(3.8,1.8)}
\pscustom[linewidth=1pt,fillstyle=solid,fillcolor=white]{
  \psellipse(7,3.5)(0.6,0.3)}
\pscustom[linewidth=1pt,fillstyle=solid,fillcolor=white]{
  \psellipse(5,3.5)(0.6,0.3)}
\pscustom[linewidth=1pt,fillstyle=solid,fillcolor=white]{
  \psellipse(9,3.5)(0.6,0.3)}

\psecurve[linewidth=5pt,linecolor=white](9,3.5)(9,1.6)(8.9,1)(8,0.5)(7,0.3)
(6,0.3)(5.5,0.4)(5,0.9)(5,1.6)(5,3.5)
\psecurve[linewidth=2pt](9,3.5)(9,1.6)(8.9,1)(8,0.5)(7,0.3)(6,0.3)(5.5,0.4)
(5,0.9)(5,1.6)(5,3.5)

\psecurve[showpoints=false,linewidth=5pt,linecolor=white](7,3)(7,3.5)(7,4)(7,4.5)(7,5)(7,5.5)(7.1,6)(7.5,6.4)
(8,6.5)(9,6.6)(10,6.6)(11,6.6)(12,6.5)(12.5,6)(12.6,5)(12.6,4)(12.6,3)(12.6,2)
(12.5,1)(12,0.5)(11,0.4)(10,0.4)(9,0.4)(8,0.4)(7.5,0.5)(7,1)(7,1.4)(7,2)
\psecurve[linewidth=2pt](7,3)(7,3.5)(7,4)(7,4.5)(7,5)(7,5.5)(7.1,6)(7.5,6.4)
(8,6.5)(9,6.6)(10,6.6)(11,6.6)(12,6.5)(12.5,6)(12.6,5)(12.6,4)(12.6,3)(12.6,2)(12.5,1)
(12,0.5)(11,0.4)(10,0.4)(9,0.4)(8,0.4)(7.5,0.5)(7,1)(7,1.4)(7,2)

\psecurve[linewidth=5pt,linecolor=white](5,2.5)(5,3.5)(5,4.5)(5,5.5)(5.2,6.5)(6,6.9)(7,7)
(8,6.9)(8.8,6.5)(9,5.5)(9,4.5)(9,3.5)(9,2.5)
\psecurve[linewidth=2pt](5,2.5)(5,3.5)(5,4.5)(5,5.5)(5.2,6.5)(6,6.9)(7,7)
(8,6.9)(8.8,6.5)(9,5.5)(9,4.5)(9,3.5)(9,2.5)

\rput(2.5,3.5){$K_1$}
\rput(4.5,6.5){$K_2$}
\rput(12,5){$K_3$}

\rput(3.6,3.5){$\beta_1$}
\rput(5,3){$\alpha_2$}
\rput(7,3){$\alpha_3$}
\rput(9,3){$\alpha_2$}

\endpspicture
\end{center}
\begin{center}
\bf Fig.~7
\end{center}

Thus $\beta_1$ maps under $H_1(T_1) \to H_1(M)$ to $0$ if and only if
all $c_j$ are even.
Otherwise it maps to an element linearly independent from the image of
$\alpha_1$.
Hence $H_1(T_1) \to H_1(M)$ is injective if and only if there is a link
component $K_j$ such that $Lk(K_1,K_j)$ is odd.
Lemma \ref{randtorus} and Theorem \ref{dichotomy} finish the proof.
\qed

The proof shows that pairs of link components with odd linking number
correspond to those pairs of ends along which the spin structure can be made
trivial.
Note that the condition on the linking numbers is extremely easy to
verify in given examples.
Since we compute modulo 2 orientations of link components are irrelevant.
If the link is given by a planar projection, then modulo 2, $Lk(K_i,K_j)$
is the same as the number of over-crossings of $K_i$ over $K_j$.

{\bf Examples.}
The complements of the following links possess a hyperbolic structure
of finite volume.
All linking numbers are even.
Hence the Dirac spectrum on those hyperbolic manifolds is discrete
for all spin structures.

\begin{center}
\pspicture(1,0)(14,10)

\psecurve[showpoints=false,linewidth=2pt](2.5,7)(2.6,6.6)(3,6.1)(3.5,5.9)
(4.1,6.2)(4.4,6.8)(4.5,7.4)(4.4,8)(4.3,8.4)
\psecurve[showpoints=false,linewidth=2pt](4.4,8)(4.3,8.4)(4.1,8.8)(3.6,9)
(3,8.8)(2.6,8.2)(2.5,7.6)(2.5,7)(2.6,6.6)

\psecurve[showpoints=false,linewidth=1pt](3.3,7.3)(3.7,7.7)(3.9,8)(4.4,8.2)
(5.1,7.8)(5.1,7.2)(4.6,6.9)(4.2,6.9)
\psecurve[showpoints=false,linewidth=1pt](4.6,6.9)(4.2,6.9)(3.8,7.2)(3.5,7.5)
(3.2,7.8)(2.7,8.1)(2.4,8.1)
\psecurve[showpoints=false,linewidth=1pt](2.7,8.1)(2.4,8.1)(1.9,7.8)(1.9,7.2)
(2.5,6.8)(3,6.9)(3.3,7.3)(3.7,7.7)

\rput(5.5,9.5){$5^2_1$}


\psecurve[showpoints=false,linewidth=2pt](9.5,7)(9.6,6.6)(10,6.1)(10.5,5.9)
(11.1,6.2)(11.4,6.8)(11.5,7.4)(11.4,8)(11.3,8.4)
\psecurve[showpoints=false,linewidth=2pt](11.4,8)(11.3,8.4)(11.1,8.8)(10.6,9)
(10,8.8)(9.6,8.2)(9.5,7.6)(9.5,7)(9.6,6.6)

\psecurve[showpoints=false,linewidth=1pt](11.3,7.1)(11.6,7)(12,7.4)(12,7.8)
(11.6,8.2)(11.1,8.1)(10.9,7.8)(10.8,7.4)
\psecurve[showpoints=false,linewidth=1pt](10.9,7.8)(10.8,7.4)(10.5,7.1)
(10.2,7.2)(10,7.8)(9.7,8)(9.4,8)
\psecurve[showpoints=false,linewidth=1pt](9.6,7.9)(9.4,8)(9,7.6)(9,7)
(9.4,6.75)(10,7.3)(10.2,7.6)
\psecurve[showpoints=false,linewidth=1pt](10,7.3)(10.2,7.6)(10.5,7.8)
(10.9,7.6)(11.3,7.1)(11.6,7)

\rput(12.5,9.5){$6^2_3$}


\psecurve[showpoints=false,linewidth=2pt](2.8,3.8)(2.9,4.05)(3.4,4.3)
(3.9,4)(4,3.4)(3.9,3)(3.8,2.8)
\psecurve[showpoints=false,linewidth=2pt](3.9,3)(3.8,2.8)(3.6,2.6)
(3.2,2.6)(2.8,3)(2.7,3.4)(2.8,3.8)(2.9,4.05)

\psecurve[showpoints=false,linewidth=1pt](4.1,3.9)(3.8,3.9)(3,3.9)
(2.4,3.8)(2,3.2)(2.4,2.4)(2.6,2.1)
\psecurve[showpoints=false,linewidth=1pt](2.4,2.4)(2.6,2.1)(3,1.6)
(3.6,1.2)(4.4,1.4)(4.5,2)(4.3,2.2)
\psecurve[showpoints=false,linewidth=1pt](4.5,2)(4.3,2.2)(3.8,3)
(3.4,3.1)(3,3)(2.8,2.7)
\psecurve[showpoints=false,linewidth=1pt](3,3)(2.8,2.7)(2.4,2.2)
(2.2,1.6)(2.8,1.1)(3.3,1.2)(3.6,1.4)
\psecurve[showpoints=false,linewidth=1pt](3.3,1.2)(3.6,1.4)(4.3,2)
(4.8,2.8)(4.8,3.4)(4.4,3.8)(4.1,3.9)(3.8,3.9)

\rput(5.5,4.5){$7^2_4$}


\psarc[linewidth=2pt](11,2){1}{70}{163}
\psarc[linewidth=2pt](11,2){1}{183}{50}

\psarc[linewidth=1pt](9.8,2){1}{65}{296}
\psarc[linewidth=1pt](9.8,2){1}{316}{45}

\psarc[linewidth=2pt,linecolor=lightgray](10.5,3){1}{296}{173}
\psarc[linewidth=2pt,linecolor=lightgray](10.5,3){1}{193}{276}

\rput(12.5,4.5){$6^3_2$}


\endpspicture
\end{center}

\begin{center}
Spectrum of the Dirac operator is discrete.
\end{center}

\begin{center}
\bf Fig.~8
\end{center}

Note that the links $5^2_1$ (Whitehead link) and $6^3_2$ (Borromeo rings)
are among the first ones for whose complements Thurston constructed
hyperbolic structures \cite{thurston97a}.

{\bf Examples.}
The complements of the following links possess a hyperbolic structure
of finite volume.
There are odd linking numbers.
Hence those hyperbolic manifolds have a spin structure for which the 
Dirac spectrum is the whole real line.

\begin{center}
\pspicture(1,0)(14,10)

\psecurve[showpoints=false,linewidth=2pt](2.6,7)(2.8,6.8)(3.4,6.4)
(3.8,6.2)(4.6,6.2)(4.7,6.8)(4.4,7)
\psecurve[showpoints=false,linewidth=2pt](4.7,6.8)(4.4,7)(4,7.2)
(3.6,7.4)(3.3,7.8)(3.4,8.2)(3.7,8.5)
\psecurve[showpoints=false,linewidth=2pt](3.4,8.2)(3.7,8.5)(3.9,8.9)
(3.2,9.3)(2.4,9)(2.1,8.2)(2.2,7.6)(2.6,7)(2.8,6.8)

\psecurve[showpoints=false,linewidth=1pt](3.5,9.1)(3.8,9.2)(5,8.6)
(5.1,7.6)(4.6,6.9)(3.8,6.4)(3.4,6.2)
\psecurve[showpoints=false,linewidth=1pt](3.8,6.4)(3.4,6.2)(2.9,6.1)
(2.6,6.2)(2.5,6.6)(2.8,7)(3.4,7.4)(3.7,7.5)
\psecurve[showpoints=false,linewidth=1pt](3.4,7.4)(3.7,7.5)(3.9,7.9)
(3.3,8.8)(3.5,9.1)(3.8,9.2)

\rput(5.5,9.5){$6^2_2$}


\psecurve[showpoints=false,linewidth=2pt](9.9,8.7)(10.3,8.7)(11,8.8)
(11.6,8.6)(12,8.2)(11.8,7.1)(11.6,6.8)
\psecurve[showpoints=false,linewidth=2pt](11.8,7.1)(11.6,6.8)(11.2,6.3)
(10.8,6)(10.1,6.1)(9.8,6.4)
\psecurve[showpoints=false,linewidth=2pt](10.1,6.1)(9.8,6.4)(9.3,7)
(9.1,7.8)(9.3,8.5)(9.9,8.7)(10.3,8.7)

\psecurve[showpoints=false,linewidth=1pt](11,6.3)(11.3,6.2)(11.9,6.6)
(11.7,7)(11.2,7.4)(10.6,7.7)(10.4,7.9)
\psecurve[showpoints=false,linewidth=1pt](10.6,7.7)(10.4,7.9)(10.1,8.7)
(10.6,9.3)(11.1,8.9)(11.1,8.6)
\psecurve[showpoints=false,linewidth=1pt](11.1,8.9)(11.1,8.6)(10.9,8.2)
(10,7.4)(9.5,7)(9.3,6.8)
\psecurve[showpoints=false,linewidth=1pt](9.5,7)(9.3,6.8)(9.2,6.4)
(9.5,6)(10.4,6.6)(11,6.3)(11.3,6.2)

\rput(12.5,9.5){$7^2_1$}


\psecurve[showpoints=false,linewidth=2pt](2.6,2)(2.8,1.8)(3.4,1.4)
(3.8,1.2)(4.6,1.2)(4.7,1.8)(4.4,2)
\psecurve[showpoints=false,linewidth=2pt](4.7,1.8)(4.4,2)(4,2.2)
(3.6,2.4)(3.3,2.8)(3.4,3.2)(3.7,3.5)
\psecurve[showpoints=false,linewidth=2pt](3.4,3.2)(3.7,3.5)(3.9,3.9)
(3.2,4.3)(2.4,4)(2.1,3.2)(2.2,2.6)(2.6,2)(2.8,1.8)

\psecurve[showpoints=false,linewidth=1pt](3.8,1.4)(3.4,1.2)(2.9,1.1)
(2.6,1.2)(2.5,1.6)(2.8,2)(3.4,2.4)(3.7,2.5)
\psecurve[showpoints=false,linewidth=1pt](3.4,2.4)(3.7,2.5)(4.65,3.2)
(4.8,3.8)(4.6,4.2)(3.9,4.2)(3.5,4.1)
\psecurve[showpoints=false,linewidth=1pt](3.9,4.2)(3.5,4.1)(3.3,4)
(3.2,3.7)(3.6,3.3)(4.4,3.1)(4.7,3)
\psecurve[showpoints=false,linewidth=1pt](4.4,3.1)(4.7,3)(4.9,2.6)
(4.8,2.2)(4.4,1.8)(3.8,1.4)(3.4,1.2)

\rput(5.5,4.5){$7^2_2$}


\psecurve[showpoints=false,linewidth=2pt](10,3)(9.7,2.9)(9.4,2.9)
(9.2,3.4)(10,3.9)(11,3.9)(11.6,3.7)(11.9,3.4)(11.9,3.1)
\psecurve[showpoints=false,linewidth=2pt](11.9,3.4)(11.9,3.1)
(11.4,2.8)(10.6,3)(10,3)(9.7,2.9)

\psecurve[showpoints=false,linewidth=1pt](9.3,3.3)(9.1,3)(9,2.4)
(9.4,1.5)(10,0.9)(10.4,0.8)(10.6,1.2)(10.5,1.5)(10.3,1.7)
\psecurve[showpoints=false,linewidth=1pt](10.5,1.5)(10.3,1.7)
(10,2.3)(9.9,2.9)(9.7,3.3)(9.3,3.3)(9.1,3)

\psecurve[showpoints=false,linewidth=2pt,linecolor=lightgray](10.3,1)(10.6,0.8)
(11.3,1)(12,2.2)(12,3.2)(11.6,3.3)(11.3,3)(11.2,2.7)
\psecurve[showpoints=false,linewidth=2pt,linecolor=lightgray](11.3,3)(11.2,2.7)
(11,2.2)(10.6,1.8)(10.2,1.4)(10.3,1)(10.6,0.8)

\rput(12.5,4.5){$6^3_1$}


\endpspicture
\end{center}

\begin{center}
For some spin structures $spec(D) = \R$.
\end{center}

\begin{center}
\bf Fig.~9
\end{center}


\providecommand{\bysame}{\leavevmode\hbox to3em{\hrulefill}\thinspace}

\ep
\ep

\parskip0ex

Universit\"at Hamburg

Fachbereich Mathematik

Bundesstr.~55

20146 Hamburg

Germany

\ep

E-Mail:
{\tt baer@math.uni-hamburg.de}

WWW:
{\tt http://www.math.uni-hamburg.de/home/baer/}

\end{document}